\documentclass{amsart}

\usepackage{amssymb}
\usepackage{amsmath}
\usepackage{graphicx}

\newtheorem{theorem}{Theorem}[section]
\newtheorem{lemma}[theorem]{Lemma}
\newtheorem{proposition}[theorem]{Proposition}
\newtheorem{proposition_definition}[theorem]{Proposition-Definition}
\newtheorem{corollary}[theorem]{Corollary}

\newtheorem{_conjecture}[theorem]{Conjecture}

\newtheorem{_problem}[theorem]{Problem}

\newtheorem{_algorithm}[theorem]{Algorithm}

\newtheorem{_subroutime}[theorem]{Subroutine}

\newtheorem{_claim}[theorem]{Claim}

\newtheorem{_subclaim}[theorem]{Sub-claim}

\newtheorem{_definition}[theorem]{Definition}
\newenvironment{definition}{\begin{_definition}\rm}{\end{_definition}}

\newtheorem{_remark}[theorem]{\it Remark}
\newenvironment{remark}{\begin{_remark}\rm}{\end{_remark}}

\newtheorem{_example}[theorem]{Example}
\newenvironment{example}{\begin{_example}\rm}{\end{_example}}

\newtheorem{_assumption}[theorem]{Assumption}
\newenvironment{assumption}{\begin{_assumption}\rm}{\end{_assumption}}

\newtheorem{_construction}[theorem]{Construction}

\numberwithin{equation}{section}
\numberwithin{table}{section}
\numberwithin{figure}{section}

\newcommand{\A}{\mathord{\mathbb A}}

\renewcommand{\P}{\mathord{\mathbb  P}}

\newcommand{\AAA}{\mathord{\mathcal A}}

\newcommand{\KKK}{\mathord{\mathcal K}}
\newcommand{\LLL}{\mathord{\mathcal L}}

\newcommand{\NNN}{\mathord{\mathcal N}}
\newcommand{\OOO}{\mathord{\mathcal O}}

\newcommand{\RRR}{\mathord{\mathcal R}}

\font\mathgot=eufm10

\newcommand{\maprightsp}[1]{\; \smash{\mathop{\; \longrightarrow \; }\limits\sp{#1}}\; }
\newcommand{\maprightsb}[1]{\; \smash{\mathop{\; \longrightarrow \; }\limits\sb{#1}}\; }

\newcommand{\mapdown}{\phantom{\Big\downarrow}\hskip -8pt \downarrow}

\newcommand{\mapdownleftright}[2]{\rlap{$\vcenter{\hbox{$\scriptstyle#1$}}$}%
\phantom{\Big\downarrow}\mapdown\rlap{$\vcenter{\hbox{$\scriptstyle#2$}}$}}
\newcommand{\mapdownsurj}{
\hbox{$\bigm\downarrow$}
\llap{\hbox{\raise 2pt\hbox{$\bigm\downarrow$}}}%
\vstrechmapdown
}

\newcommand{\mapup}{\phantom{\Big\uparrow}\hskip -8pt \uparrow}
\newcommand{\mapupright}[1]{\mapup\rlap{$\vcenter{\hbox{$\scriptstyle#1$}}$}}

\newcommand{\mapupsurj}{
\hbox{$\bigm\uparrow$}
\llap{\hbox{\raise 2pt\hbox{$\bigm\uparrow$}}}%
\vstrechmapup
}

\newcommand{\inj}{\hookrightarrow}

\newcommand{\isom}{\smash{\mathop{\;\to\;}\limits\sp{\sim\,}}}

\newcommand{\set}[2]{\{\; {#1} \; \mid \; {#2} \;  \}}

\newcommand{\map}[3]{ #1 \, : \, #2 \, \to \, #3}
\newcommand{\mapisom}[3]{ #1 \, : \, #2 \; \isom \; #3}

\newcommand{\shortmap}[3]{ #1 : #2 \to #3}

\newcommand{\sm}{\setminus}
\newcommand{\st}{\subset}

\newcommand{\wt}{\widetilde}

\newcommand{\sprime}{\sp\prime}

\newcommand{\spcirc}{\sp{\mathord{\circ}}}

\newcommand{\dual}{\sp{\vee}}
\newcommand{\comp}{\sp{\wedge}}

\newcommand{\inv}{\sp{-1}}

\newcommand{\Ker}{\operatorname{\rm Ker}\nolimits}

\newcommand{\pr}{\operatorname{\rm pr}\nolimits}

\newcommand{\GL}{\operatorname{\it GL}\nolimits}

\newcommand{\SL}{\operatorname{\it SL}\nolimits}

\newcommand{\Hom}{\operatorname{\rm Hom}\nolimits}

\newcommand{\sign}{\operatorname{\rm sign}\nolimits}
\newcommand{\Sing}{\operatorname{\rm Sing}\nolimits}
\newcommand{\Spec}{\operatorname{\rm Spec}\nolimits}

\newcommand{\rank}{\operatorname{\rm rank}\nolimits}

\newcommand{\ch}{\operatorname{\rm char}\nolimits}

\newcommand{\closure}[1]{\overline{#1}}

\newcommand{\Der}[2]{\frac{\partial #1}{\partial #2}}
\newcommand{\dDer}[2]{\displaystyle\Der{#1}{#2}}
\newcommand{\DDer}[3]{\frac{\partial^2 #1}{\partial #2 \partial #3}}
\newcommand{\dDDer}[3]{\displaystyle\DDer{#1}{#2}{#3}}

\newcommand{\rmand}{\textrm{and}}

\newcommand{\qquand}{\qquad\rmand\qquad}
\newcommand{\quand}{\quad\rmand\quad}

\newcommand{\mystruth}[1]{\phantom{\hbox{\vrule height #1}}}
\newcommand{\mystrutd}[1]{\phantom{\hbox{\vrule depth #1}}}

\newcommand{\dsum}{\displaystyle\sum}

\newcommand{\XXX}{\mathord{\rlap{\hbox{\,$\overline{\hbox to 7pt{\phantom{I}}}$}}X}}
\newcommand{\XX}{\mathord{X}}
\newcommand{\X}{\mathord{X\spcirc}}
\newcommand{\DDD}{\mathord{\closure{D}}}
\newcommand{\DD}{\mathord{D}}

\newcommand{\fDD}{\mathord{\mathcal D}}
\newcommand{\fD}{\mathord{\mathcal D\spcirc}}
\newcommand{\fCC}{\mathord{\mathcal C}}
\newcommand{\fC}{\mathord{\mathcal C\spcirc}}
\newcommand{\fE}{\mathord{\mathcal E}}
\newcommand{\fR}{\mathord{\mathcal R}}
\newcommand{\fH}{\mathord{\mathcal H}}
\newcommand{\fHt}{\fH\,\tilde{}}

\newcommand{\Degasymm}{\mathord{\mathbf R}}

\newcommand{\Tan}[1]{\mathord{T}(#1)}

\newcommand{\PM}{\mathord{{\mathbf P}}}

\newcommand{\XXP}{\XX\times \PM}
\newcommand{\XP}{\X\times \PM}
\newcommand{\pf}{(p, [f])}
\newcommand{\DDDf}{\DDD_{[f]}}
\newcommand{\DDf}{\DD_{[f]}}

\newcommand{\Bs}{\operatorname{\rm Bs}\nolimits}
\newcommand{\rest}{\,|\,}

\font\bfti=cmbxti10
\newcommand{\vb}{\mathord{\hbox{\bfti b}}}
\newcommand{\tvb}{\tilde{\vb}}
\newcommand{\fRsm}{\fR\sp{\mathord{\mathrm{sm}}}}

\newcommand{\vpism}{\varpi\sp{\mathord{\mathrm{sm}}}}

\newcommand{\mmm}{\mathord{\hbox{\mathgot m}}}

\renewcommand{\SL}{S_L}
\newcommand{\CL}{C_\Lambda}
\newcommand{\wCL}{\widetilde{C}_\Lambda}

\newcommand{\AAeven}{\mathord{(\widetilde{\text{\rm A}})}}       
\newcommand{\typeC}{\mathord{(\text{\rm C})}}
\newcommand{\typeR}{\mathord{(\text{\rm R})}}

\newcommand{\threevars}{ x_1, x_2, z}
\newcommand{\sixvars}{ u_1, v_1, u_2, v_2, w, y}
\newcommand{\dimV}{l}

\title[Dual  varieties  in characteristic $2$]{
Singularities of  dual  varieties in characteristic $2$}

\author[I.~Shimada]{Ichiro Shimada}
\address{
Department of Mathematics\\
Faculty of Science\\
Hokkaido University\\
Sapporo 060-0810\\
JAPAN
}
\email{shimada@math.sci.hokudai.ac.jp
}


\thanks{Partially supported by *****}

\begin{document}

\begin{abstract}
We investigate  unibranched singularities of  dual varieties of
  even-dimensional smooth projective varieties in characteristic $2$.
\end{abstract}

\maketitle
\section{Introduction}\label{sec:intro}

We work over an algebraically closed field $k$.
\par
\medskip
Let $X\subset \P^m$ be a smooth projective variety of dimension $n>0$.
We assume that $X$ is not contained in any hyperplane of $\P^m$.
Then we can consider
 the dual projective space 
$$
\PM:=(\P^m)\dual
$$
of $\P^m$  as the parameter space of the linear system $|M|$ of hyperplane sections of $X$,
where $M$ is a linear subspace of $H^0 (X, \OOO_X(1))$.
We use the same letter to denote a point $H\in \PM$ and the corresponding hyperplane $H\st \P^m$.
Let $\fDD\subset X\times \PM$ be the universal family of the hyperplane sections.
Then $\fDD$ is a smooth scheme of dimension $n+m-1$. 
The support of $\fDD$ is equal to the set
$$
\set{(p, H)\in X\times  \PM}{p\in H}.
$$
Let $\fCC\subset \fDD$ be the critical subscheme of the second projection
$\fDD\to \PM$.
(See  Notation and Terminology below for the definition of critical subschemes.)
Then $\fCC$ is smooth, irreducible and  of dimension $m-1$.
If $\NNN$ is the conormal sheaf of $X\subset \P^m$,
then $\fCC$ is isomorphic  to $\P ^* (\NNN)$ (\cite[Remarque 3.1.5]{MR0354657Katz}).
The support of $\fCC$ is equal to the set
$$
\set{(p, H)\in \fDD}{\hbox{the divisor $H\cap X$ of $X$  is singular at $p$}}.
$$
We denote the projections by
$$
\pi_1: \fCC\to X\quand \pi_2: \fCC\to\PM.
$$
The image of $\pi_2$ is called 
the \emph{dual variety} of $X\subset \P^m$.
Let $\fE\subset \fCC$ be the critical subscheme of 
$\pi_2$.
By~\cite[Proposition 3.3]{MR0354657Katz}, 
 $\fE$ is set-theoretically equal to
$$
\set{(p, H)\in \fCC}{\hbox{the Hessian of the  singularity  of  $H\cap X$  at $p$ is degenerate}}.
$$
We will study the singularities of the dual variety by investigating the morphism $\pi_2$
at points of $\fE$.
\par
\medskip
Since the paper of Wallace~\cite{MR0080354},
properties   of dual varieties peculiar to  positive characteristics
have been studied mainly from the point of
view of the (failure of the) reflexivity.
See Kleiman's paper~\cite{MR846021} for the definition and a  detailed account of the reflexivity.
Many studies have been done for the analysis of the situation
in which the reflexivity does not hold.
See~\cite{ MR882943, MR986811,   MR895153, MR1143223,  MR996326, MR1998935, MR1427663,   MR867952}, for example.
\par
\medskip
A well-known  example of  the situations in which the reflexivity fails is as follows.
Suppose that
\begin{equation}\label{eq:char2odd}
\ch k=2
\quand 
\dim X\equiv 1 \bmod 2.
\end{equation}
Then the critical subscheme  $\fE$ of $\pi_2:\fCC\to\PM$ coincides with $\fCC$
(\cite[Section 1.2]{MR0354657Katz}, \cite[Corollary (18)]{MR846021}),
and hence either the dual variety is not a hypersurface in $\PM$,  or $X$ is not reflexive
by the  generalized Monge-Segre-Wallace criterion
(\cite[Theorem (4.4)]{MR775882}, \cite[Theorem (4)]{MR846021}).
\par
\medskip
Except for the case~\eqref{eq:char2odd}, however,
the reflexivity is recovered when the linear system $|M|$ of the hyperplane sections are sufficiently ample.
We have the following theorem (\cite[Th\'eor\`eme 2.5]{MR0354657Katz}, \cite[Theorem (5.4)]{MR895153}):
\begin{theorem}
Suppose that $\ch k\ne 2$ or $\dim X$ is even.
If $X$ is embedded in $\P^m$ by a complete linear system of the form 
 $|\AAA\sp{\otimes d}|$ with $\AAA$ being  a very ample line bundle and $d\ge 2$,
then the dual variety of $X\subset \P^m$ is a hypersurface of $\PM$,
and $X\subset \P^m$ is reflexive.
\end{theorem}
See also~\cite[Proposition 4.9]{char3} or Proposition~\ref{prop:eval} of this paper
for other sufficient conditions 
for the dual variety to be  a hypersurface, and 
for the reflexivity to hold.
\par
\medskip
In our previous paper~\cite{char3},
we have discovered that,
even when the reflexivity holds,
the singularities of dual varieties in characteristic $3$
still possess a peculiar feature.
We assume that the linear system $|M|$  is sufficiently ample.
In particular,
the dual variety is a hypersurface in $\PM$ and $X$ is reflexive.
We consider the projective plane curve obtained by cutting the dual variety by a general plane in $\PM$.
If $k$ is of characteristic $>3$ or $0$,
then the plane curve has only ordinary cusps as its unibranched singularities.
In characteristic $3$, 
the plane curve has  singular points of
$E_6$-type instead   of ordinary cusps.
\par
\medskip
In this paper, we investigate the singularities of  the dual  variety 
in the case where $\ch k=2$ and $\dim X$ is even.
\par
\medskip
The dual variety is the \emph{discriminant variety} 
associated with the linear system of hyperplane sections.
See~\cite{MR644816} for the definition of discriminant varieties.
In fact,
our results are proved 
not only for dual varieties 
but for discriminant varieties
associated with (not necessarily very ample) linear systems.
Here in Introduction,
however,
we present our results 
for dual varieties.
\par
\medskip
Let $k$ be of characteristic $2$, and let $\dim X$ be even.
For simplicity,
we assume that $|M|$ is  sufficiently ample so that 
 the evaluation homomorphism
$$
\map{v_p\sp{[4]}}{M}{\LLL_p/\mmm_p^5 \LLL_p}
$$
is surjective at every point $p$ of $X$,
where $\mmm_p$ is the maximal ideal of   $\OOO_{X, p}$,
and $\LLL_p$ is the $\OOO_{X,p}$-module $\LLL\otimes \OOO_{X, p}$.
Under this assumption, the critical subscheme 
 $\fE$ is an irreducible 
divisor of $\fCC$
and hence the dual variety is a hypersurface of $\PM$.
\par
\medskip
We can define the \emph{universal Hessian}
$$
\map{\fH}{\pi_1 ^* \Tan{X} \otimes \pi_1 ^* \Tan{X} }{ \pi_1^* \LLL \otimes \pi_2^* \OOO_{\PM} (1) }
$$
 on $\fCC$,
where $\Tan{X}$ is the tangent bundle of $X$.
As was proved in~\cite[Proposition 3.14]{char3},
the  subscheme $\fE$ coincides with the degeneracy subscheme of the homomorphism
$$
\pi_1 ^* \Tan{X}\;\;\to\;\; \pi_1^* \LLL \otimes \pi_2^* \OOO_{\PM} (1)\otimes (\pi_1 ^* \Tan{X})\dual
$$
induced from  $\fH$,
where $(\pi_1 ^* \Tan{X})\dual$ is the dual vector bundle of $\pi_1 ^* \Tan{X}$.
(See Notation and Terminology below for the definition of degeneracy subschemes.)
The peculiarity of the geometry of the dual variety in characteristic $2$
stems from the fact that the universal Hessian is 
not only symmetric but also \emph{anti}-symmetric;
that is, $\fH(x\otimes x)=0$ and $\fH(x\otimes y)+\fH(y\otimes x)=0$ hold for any local sections $x$ and $y$ of $\pi_1 ^* \Tan{X}$.
From this fact, it follows that the irreducible divisor 
$\fE$   is written as $2\,\fR$,
where $\fR$ is a reduced divisor of $\fCC$.
We denote by $\varpi_2 :\fR\to\PM$ the projection.
\begin{definition}
We put $R:=k[[t_1, \dots, t_{m-4}]]$.
It turns out that the equations
\begin{multline}\label{eq:defeqGamma}
u_1 y^2+v_1\;\;=\;\;u_2 y^2+v_2\;\;=\;\;w v_1 + u_1 v_2 y \;\;=\;\; w v_2 + u_2 v_1 y \;\;=\;\;\\
\;\;=\;\; v_1 v_2 + w y \;\;=\;\;w^2+u_1 u_2 y^2 \;\;=\;\;0\phantom{aaaa}
\end{multline}
define an $(m-1)$-dimensional singular scheme in 
the $(m+2)$-dimensional smooth scheme $\Spec R[[ u_1, v_1, u_2, v_2, w, y]]$,
which we will denote by $\Gamma_{m-1}$.
\end{definition}
The scheme  $\Gamma_{m-1}$ is singular  along the  $(m-2)$-dimensional  locus defined by
$$
v_1=v_2=w=y=0.
$$
See~\S\ref{sec:Gamma} for the geometric meaning of $\Gamma_{m-1}$.
\par
\medskip
Let $P=(p, H)\in \XXP$ be a general point of the reduced irreducible divisor $\fR$ of $\fCC$.
Then we have the following:
\begin{itemize}
\item[(I)]
There exist isomorphisms of local rings
$$
\OOO_{\pi_2\inv (H), P}\cong k[[s, t]]/(s^2, t^2)
\qquand
\OOO_{\varpi_2\inv (H), P}\cong k[[s, t]]/(s^2, t^2)
$$ 
over $k$.
In particular,
the projection $\varpi_2 :\fR\to\PM$ is inseparable of degree $4$ over its image.
\item[(II)]
The formal completion 
$$
\map{(\pi_2)_P\comp}{\Spec (\OOO_{\fCC, P})\comp}{\Spec (\OOO_{\PM, H})\comp}
$$
 of the projection  $\pi_2 :\fCC\to\PM$ at $P$
factors through a singular scheme  isomorphic to $\Gamma_{m-1}$.
\item[(III)]
Let $\Lambda$ and $L$ be general linear subspaces of $\PM$ of dimension $2$ and $3$, respectively, 
such that $H\in \Lambda\st L$.
We put 
$\CL:=\pi_2\inv (\Lambda)$ and $\SL:=\pi_2\inv (L)$.
Then $\SL$ is smooth of dimension $2$ at $P$,
and $\CL$ is a curve on $\SL$ that has an ordinary cusp at $P$.
\item[(IV)]
Let $\nu : \wCL\to \CL$ be the normalization of $\CL$ at $P$,
and let $z$  be a formal parameter of $\wCL$ at the inverse image $P\sprime\in\wCL$ of $P$.
Then the formal completion at $P\sprime$ of the composite
of $\nu : \wCL\to \CL$  and  $\pi_2\rest \CL: \CL\to\Lambda$
is written as
$$
\renewcommand{\arraystretch}{1.28}
\begin{array}{lclcl}
((\pi_2\rest\CL)\circ\nu)^* x &=& a\, z^4 &+&  (\hbox{terms of degree $\ge 6$})
\quand \\
((\pi_2\rest\CL)\circ\nu)^* y &=& b\, z^4 &+&  (\hbox{terms of degree $\ge 6$})
\end{array}
$$
for some $a, b\in k$,
where $(x, y)$ is a formal parameter system of $\Lambda$ at $H$.
\end{itemize}
\par
\medskip
This paper is organized as follows.
First we define several notions in~\S\ref{sec:preliminaries}.
The conditions on $P=(p, H)\in \fR$ for the facts (I)-(IV) above to hold will be stated more precisely
in terms of the singularity of the divisor $H\cap X$ of $X$ at $p$.
For this purpose, 
we define some classes of hypersurface singularities  in  \S\ref{sec:HS}.
In \S\ref{sec:Gamma},
the $(m-1)$-dimensional singular scheme  $\Gamma_{m-1}$ is introduced.
After recalling  the definitions  and  results given in \cite[Section 3]{char3} in \S\ref{sec:reminder}, 
we prove the main results above in more refined forms in~\S\ref{sec:evenmain}.
In~\S\ref{sec:div}, we give a remark about the degree of $\fR$
with respect to $\OOO_{\PM}(1)$,
and derive a nontrivial divisibility relation 
among Chern numbers of $X$ from the fact (I) above.
\par
\medskip
The author would like to thank Professor Hajime Kaji for many valuable comments
and suggestions.
\par
\bigskip
{\bf Notation and Terminology.}
\begin{itemize}
\item[(1)] We work over an algebraically closed field $k$.
By a \emph{variety}, we mean  a reduced irreducible quasi-projective scheme over $k$. 
A \emph{point} of a variety means a closed point.
Let $X$ be a  variety, and $P \in X$ a  point.
We  denote by $(X, P)\comp$ the scheme $\Spec (\OOO_{X, P})\comp$,
where $(\OOO_{X, P})\comp$ is the formal completion of the local ring $\OOO_{X, P}$ of $X$ at $P$.
\item[(2)] 
Let $X$ be a smooth variety.
We denote by $T_P (X)$ the Zariski tangent space to $X$ at a point $P\in X$,
and  by $\Tan{X}$ the tangent bundle of $X$.
\item[(3)]
Let $E$ and $F$ be vector bundles on a variety $X$ with rank $e$ and $f$,
respectively,
and let $\sigma : E\to F$ be a bundle homomorphism.
We put $r:=\min (e, f)$.
The \emph{degeneracy subscheme} of $\sigma$
is  the closed subscheme of $X$ defined locally on $X$ by all the $r$-minors 
of the $f\times e$-matrix expressing $\sigma$.
\item[(4)]
Let $f: X\to Y$ be a morphism from a smooth variety $X$ to a smooth variety $Y$.
The  \emph{critical subscheme} of $f$ is
the degeneracy subscheme of the  homomorphism
$\shortmap{df}{\Tan{X}}{f\sp *\,\Tan{Y}}$.
\item[(5)]
For a formal power series $F$ with coefficients in $k$,
we denote by $F^{[d]}$ the homogeneous part of degree $d$ in $F$.
\end{itemize}
%
%
%\newpage
%
%
\section{Preliminaries}\label{sec:preliminaries}
\subsection{The quotient morphism by an integrable tangent subbundle}\label{subsec:quot}
Let $X$ be a smooth variety defined over an algebraically closed field of characteristic $p>0$.
\begin{definition}
A subbundle $\NNN$ of $\Tan{X}$ is called \emph{integrable}
if $\NNN$ is closed under the  $p$-th power operation 
and the bracket product of derivations.
\end{definition}
\begin{proposition}[\cite{MR0450263Exp6} Th\'eor\`eme 2]
Let  $\NNN$ be an integrable  subbundle of $\Tan{X}$.
Then there exists  a unique morphism
$q: X\to X\sp{\NNN}$  with the following properties;
\begin{itemize}
\item[(i)] $q$ induces a homeomorphism on the underlying topological spaces, 
\item[(ii)] $q$ is a radical covering of height $1$, and 
\item[(iii)] the kernel of 
$dq : \Tan{X}\to q^* \,\Tan{X\sp{\NNN}}$ coincides with $\NNN$.
\end{itemize}
Moreover, the variety $X^{\NNN}$ is smooth,  and the morphism $q$ is finite of degree $p^r$, where $r$ is the rank of $\NNN$.
\end{proposition}
\begin{definition}
For an integrable subbundle $\NNN$ of $\Tan{X}$,
the morphism $q: X\to X^{\NNN}$ is called 
the  \emph{quotient morphism}
 by $\NNN$.
\end{definition}
From the construction of the quotient morphism given in the proof of \cite[Th\'eor\`eme 2]{MR0450263Exp6},
we obtain the following:
\begin{proposition}\label{rem:facquot}
Let  $f: X\to Y$ be a  morphism from a smooth variety $X$ to a smooth variety $Y$
such that 
the kernel
$\KKK$ of  $df : \Tan{X}\to  f^*\,\Tan{Y}$ 
is a subbundle of $\Tan{X}$.
Then $\KKK$ is integrable,
and the morphism $f$ factors through the quotient morphism $q: X\to X\sp{\KKK}$ by $\KKK$.
\end{proposition}
\subsection{Anti-symmetric forms}\label{subsec:anti-symmetric}
\begin{definition}\label{def:anti-symmetric_form}
Let $V$ be  a finite dimensional vector space over a field $K$.
A bilinear form $\varphi : V\times V\to K$
is called \emph{anti-symmetric} if the following conditions are satisfied:
\begin{itemize}
\item[(i)] $\varphi (v, v)=0$ for any $v\in V$, and 
\item[(ii)] $\varphi (v, w)+\varphi(w, v)=0$ for any $v, w\in V$.
\end{itemize}
Note that, in characteristic $2$, the condition (ii) does not imply the condition (i).
\end{definition}
\begin{definition}\label{def:anti-symmetric_matrix}
An $n\times n$-matrix $A=(a_{ij})$ with components in a commutative ring is called \emph{anti-symmetric}
if the following conditions are satisfied:
\begin{itemize}
\item[(i)] $a_{ii}=0$ for $i=1, \dots, n$, and 
\item[(ii)] $a_{ij}+a_{ji}=0$ for $i, j=1, \dots, n$.
\end{itemize}
An  anti-symmetric matrix in characteristic $2$ is just a symmetric matrix with zero diagonal components.
\end{definition}
The following results are  well-known.
We put
$$
J_2:=\left[
\begin{array}{cc}
0 & 1 \\
-1 & 0 
\end{array}
\right],
$$ 
and let 
$J_{2r, n}$ be the $n\times n$ matrix
$$
J_{2r, n}:=\left[
\begin{array}{cccc}
J_2 & & & \\
 &\smash{\ddots} & &\smash{\raise 5pt \hbox{\huge $0$}} \\
 & & \smash{J_2}& \\
 & \hskip -20pt \hbox{\huge $0$}& &\hbox{\huge $0$}
\end{array}
\right]
$$
with $r$ copies of $J_2$ along the diagonal.
\begin{lemma}\label{lem:J2}\label{lem:anti-symmetric_rank}
Let $A$ be an anti-symmetric matrix of type $n\times n$
with  components  in a field $K$.
Then the rank of $A$ is even, and there exists $T\in \GL(n, K)$ such that
${}\sp t  T A T$ is equal to $J_{2r, n}$,
where $2r$ is the rank of $A$.        
\end{lemma}
For a positive integer $m$,
we define a homogeneous polynomial $f_{2m}$ of degree $m$ in variables $x_{ij}\; (1\le i<j\le 2m)$ by
$$
f_{2m}:=\sum\,
\sign\left(
\begin{array}{ccccc}1& 2 &\dots &&2m \\ i_1 &j_1 &\dots &i_m &j_m\end{array}
\right)
x_{i_1 j_1} x_{i_2 j_2} \cdots x_{i_mj_m},
$$
where the summation ranges over  all the lists
$[\, [i_1, j_1], \dots, [i_m, j_m]\,]$
satisfying the conditions $i_\nu<j_\nu \; (\nu=1, \dots, m)$, $i_1<i_2<\dots<i_m$,  and
$$
\{ i_1, j_1, \dots, i_m, j_m\}=\{1, 2, \dots, 2m\}.
$$
\begin{lemma}\label{lem:square}
Let 
\begin{equation}\label{eq:Ax}
A=
\left[
\begin{array}{ccccc}
0 & x_{12} & x_{13} &  \;\;\;\dots &\\
-x_{12} & 0 & x_{23}  & \;\;\;\dots  &\\
-x_{13} & - x_{23} & 0 &  \;\;\;\dots & \\
&\cdots& &\;\;\;\ddots &
\end{array}
\right]
\end{equation}
be an anti-symmetric matrix of type $2m\times 2m$
with  components  being variables  $x_{ij} \;(i<j)$.
Then $\det A$ is equal to $f_{2m}^2$.
\end{lemma}
For a positive integer $i$, we put
$$
\tau(i):=\begin{cases}
i-1 & \textrm{if $i$ is even,}\\
i+1 & \textrm{if $i$ is odd.}
\end{cases}
$$
\begin{corollary}\label{cor:sqrtdetH}
Let $A$ be an anti-symmetric matrix as in Lemma~\ref{lem:square}.
We put
$$
g_A:=\sqrt{\;\det(J_{2r, 2m}+A)\;},
$$
which is a polynomial  of $x_{ij}$ by Lemma~\ref{lem:square}.
\par
{\rm (1)}
If $2r\le 2m-4$,
then 
$g_A$ has no terms of degree $\le 1$.
\par
{\rm (2)}
Suppose that $2r=2m-2$.
Then $g_A$ is equal, up to sign,  to 
$$
z  
\;+\; \sum_{j=1}^{m-1} \varepsilon_j\,x_{2j-1, 2j}\,z 
\;+\;  \sum_{i=1}^{2m-2} \varepsilon\sprime_i\,x_{i, 2m-1}\,x_{\tau(i), 2m}
\;+\;(\hbox{\rm terms of degree $\ge 3$}),
$$
where $z:=x_{2m-1, 2m}$, and  $\varepsilon_j$ and $\varepsilon\sprime_i$ are $\pm 1$.
\end{corollary}
Using~\eqref{eq:Ax}, we consider the affine space 
$$
\A^{m(2m-1)}=\Spec k[\dots,  x_{ij}, \dots ]\;\;\;(1\le i<j\le 2m)
$$
as the space of anti-symmetric matrices of type $2m\times 2m$.
\begin{lemma}\label{lem:irredfm}
The hypersurface in $\A^{m(2m-1)}$ defined by $f_{2m}=0$ is irreducible.
\end{lemma}
\begin{proof}
By Lemma~\ref{lem:J2},
the hypersurface defined by $f_{2m}=0$  is the closure of the locus
$$
\set{{}^t T J_{2m-2,2m} T}{T\in \GL(2m, k)},
$$
and hence is irreducible.
\end{proof}
\subsection{The formula of Harris-Tu-Pragacz}\label{subsec:HTP}
Let $X$ be a smooth variety,
$E$ a vector bundle on $X$,
and $L$ a line bundle on $X$.
\begin{definition}\label{def:anti-symmetricE}
A bundle homomorphism
$$
\map{\sigma}{E\otimes \sb{\OOO_{X}} E}{L}
$$
is called \emph{anti-symmetric} 
if the following conditions are satisfied:
\begin{itemize}
\item[(i)] $\sigma (x\otimes x)=0$ for any local section $x$ of $E$, and 
\item[(ii)] $\sigma (x\otimes y)+\sigma (y\otimes x)=0$ for any local sections $x$ and $y$ of $E$.
\end{itemize}
\end{definition}
\begin{definition}\label{def:degeneracy_anti-symmetric}
For an anti-symmetric homomorphism 
$\shortmap{\sigma}{E\otimes E}{L}$,
we define the \emph{degeneracy subscheme} $\Degasymm(\sigma)$  as follows.
Let $U$ be an arbitrary Zariski open subset of $X$ over which $E$ and $L$ are trivialized.
Let $A_U$ be the anti-symmetric matrix with  components  in $\Gamma (U, \OOO_X)$
expressing $\sigma$ over $U$.
By Lemma~\ref{lem:square}, there exists $f_U\in \Gamma (U, \OOO_X)$
such that $\det A_U=f_U^2$.
Then   the collection  of the closed subschemes $\{f_U=0\}$ of Zariski open subsets $U$
patches together to form  a closed  subscheme $\Degasymm(\sigma)$ of $X$,
which is the degeneracy subscheme of $\sigma$.
\end{definition}
\begin{definition}
For a closed subscheme $W$ of $X$,
we denote by
$[W]$
the class of $W$ in the Chow group  of $X$.
\end{definition}
The following was proved by Harris-Tu~\cite{MR721453} 
in characteristic $0$ and by Pragacz~\cite{MR1171269} in any characteristics.
\begin{theorem}\label{thm:HTP}
Suppose that the rank $e$ of $E$ is even, and that $\Degasymm(\sigma)$ is of codimension $1$ in $X$.
Then we have
$$
[\Degasymm(\sigma)]\;=\;\left(\,\frac{e}{2}\,c_1 (L) \,-\, c_1 (E)\,\right)\;\cap\; [X]
$$
in the Chow group of $X$.
\end{theorem}
\begin{remark}\label{rem:as_a_divisor}
Suppose that $\Degasymm(\sigma)$ is of codimension $1$ in $X$.
Let
$$
\map{\sigma\,\tilde{}}{E}{L\otimes E\dual}
$$
be the homomorphism induced from  $\sigma$,
where $E\dual$ is the dual vector bundle of $E$.
Then the degeneracy subscheme 
of  $\sigma\,\tilde{}\,$ is also a divisor of $X$,
and is equal to $2\,\Degasymm(\sigma)$ 
as a divisor of $X$.
\end{remark}
\subsection{Ordinary cusps}\label{subsec:cusp}
Let $P$  be a point of a smooth surface $S$,
and let $(s, t)$ be a formal parameter system of $S$ at $P$.
For $\phi\in  (\OOO_{S, P})\comp= k[[s, t]]$,
we denote by $\phi^{[d]}$  the homogeneous part of degree $d$ in $(s, t)$.
Let $C$ be a divisor of $S$ that contains $P$ and is singular at $P$.
\begin{definition}\label{def:degenerate}
We say that $\phi\in k[[s, t]]$ is \emph{degenerate}
if $\phi^{[0]}=\phi^{[1]}=0$ and $\phi^{[2]}=l(s, t)^2$ hold for some   linear form $l(s, t)$ of $(s, t)$.
It is obvious that this definition of degeneracy  does not depend on the choice of
the   formal parameter system $(s, t)$.
We say that $C$ has a \emph{degenerate singularity} at $P$
if a (and hence any) formal power series defining $C$ at $P$ is degenerate.
\end{definition}
\par
\medskip
The equivalence of the following conditions on the singularity of $C$ at $P$,
which is well-known in characteristic $0$,
 holds also in characteristic $2$.
\begin{itemize}
\item[(i)]
For an arbitrary formal parameter system $(s, t)$ of $S$ at $P$,
 any formal power series $\phi$ defining  $C$ at $P$ is degenerate
and the linear form $\sqrt{\phi^{[2]}}$ is not zero and does not divide $\phi^{[3]}$.
\item[(ii)]
There exist a formal parameter system $(s, t)$ of $S$ at $P$
and a defining formal power series $\phi$ of $C$ at $P$
such that
 $\phi$  is degenerate and that
 the linear form $\sqrt{\phi^{[2]}}$ is not zero and does not divide $\phi^{[3]}$.
\item[(iii)]
There exists a formal parameter system $(s, t)$ of $S$ at $P$
such that $C$ is defined  by $s^2+t^3=0$ locally at $P$.
\end{itemize}
\begin{definition}
We say that $C$ has an \emph{ordinary cusp} at $P$ if the conditions (i)-(iii) above are satisfied.
\end{definition}
\section{Hypersurface singularities in characteristic $2$}\label{sec:HS}
From now on until the end of the paper,
we assume that the base field $k$ is of characteristic $2$.
\par
\medskip
Let $X$ be a smooth variety of dimension $n$.
Let $p$ be a point of $X$,
and let $D$ be a hypersurface of $X$ that passes through $p$ and is singular at $p$.
Let $(x_1, \dots, x_n)$ be a formal parameter system of $X$ at $p$.
Suppose that $D$ is defined by $\phi=0$ locally at $p$,
where $\phi$ is a formal power series of $(x_1, \dots, x_n)$.
The $n\times n$ matrix
$$
H_{\phi, p}:=\left(\DDer{\phi}{x_i}{x_j}  (p)\right)
$$
defining the Hessian
$$
T_p (X)\times T_p(X)\;\;\to\;\; k
$$
of $D$ at $p$ is \emph{anti}-symmetric because we are in characteristic $2$.
\begin{definition}
A formal parameter system $(x_1, \dots, x_n)$ of $X$ at $p$ is said to be \emph{admissible with respect to $\phi$} if
the matrix $H_{\phi, p}$ is of the form $J_{2r, n}$, where $2r$ is the rank of the Hessian of $D$ at $p$.
\end{definition}
\begin{remark}\label{rem:make_admissible}
By Lemma~\ref{lem:J2},
a formal parameter system admissible with respect to $\phi$
is obtained from an arbitrarily given formal parameter system
by a \emph{linear} transformation of parameters.
\end{remark}
From now on to the end of this section,
we assume that $n$ is even.
\par
\medskip
We put
$$
\rho:=\;\;f_{n} \left(\dots,\;\;\DDer{\phi}{x_i}{x_j},\;\; \dots \right)\;\;=\;\;\sqrt{\det\left(\DDer{\phi}{x_i}{x_j}\right)}\;\;,
$$
where $f_n$ is the polynomial defined in~\S\ref{subsec:anti-symmetric}.
\begin{definition}\label{def:CR}
Let $C(D, p)$ be the subscheme of $(X, p)\comp=\Spec (\OOO_{X, p})\comp$ defined by
$$
\phi=\Der{\phi}{x_1}=\dots=\Der{\phi}{x_n}=0,
$$
and let $R(D, p)$  be the subscheme of $(X, p)\comp$ defined by
$$
\phi=\Der{\phi}{x_1}=\dots=\Der{\phi}{x_n}=\rho=0.
$$
It is obvious that each of $C(D, p)$ and $R(D, p)$ depends only on $D$ and $p$,
and not on the choice of $\phi$ and $(x_1, \dots, x_n)$.
\end{definition}
Using the admissible formal parameter system with respect to $\phi$,
or the subschemes $C(D, p)$ and $R(D, p)$ of $(X, p)\comp$,
we can define various classes of hypersurface singularities in characteristic $2$.
\par
\medskip
Suppose  that the rank of the Hessian matrix $H_{\phi, p}$ is $n-2$.
Let
$$
(x, \xi)=(x_1, \dots, x_{n-2}, \xi_1, \xi_2)
$$
be an admissible formal parameter system with respect to $\phi$.
We put
$$
G:=\set{M\in \GL(n-2, k)}{{}^t M J_{n-2, n-2} M=J_{n-2, n-2}}.
$$
The following lemma is trivial:
\begin{lemma}\label{lem:trans}
Let 
$$
(x\sprime, \xi\sprime)=(x_1\sprime, \dots, x_{n-2}\sprime, \xi\sprime_1, \xi\sprime_2)
$$
be another formal parameter system of $X$ at $p$.
Then $(x\sprime, \xi\sprime)$ is admissible with respect to $\phi$ if and only if 
$(x, \xi)$ and $(x\sprime, \xi\sprime)$ are related by the transformation of the form
$$
\left[
\begin{array}{c}
\\
x \\
\\
\hline
\xi\mystruth{12pt}
\end{array}
\right]
=\left[
\begin{array}{ccc|c}
&& & \\
&A& & 0 \\
&& & \\
\hline 
&*& & B \mystruth{12pt}
\end{array}
\right]
\left[
\begin{array}{c}
\\
x \sprime\\
\\
\hline
\xi\sprime\mystruth{12pt}
\end{array}
\right]
+
(\hbox{\rm terms of degree $\ge 2$})
$$
such that $A\in G$ and $B\in \GL(2, k)$.
\end{lemma}
We denote the coefficients in the formal power series expansion  of $\phi$
in the parameters $(x, \xi)$ as indicated in Table~\ref{table:coeff}.
\begin{table}
$$
\renewcommand{\arraystretch}{1.4}
\begin{array}{lllllllllll}
\phantom{+}&\rlap{$x_1\, x_2 \,+\, \dots \,+\,x_{n-3}\, x_{n-2}\,+$} &&&&&&&& &\\
+& a_1 \, x_1 ^2 &+& \phantom{\cdots}\cdots\phantom{\cdots}& + &a_{n-2}\, x_{n-2}\,^2 & +&\alpha_1 \xi_1^2 &+&\alpha_2\xi_2^2&+\\
+& b_1 \,x_1\, \xi_1^2& +& \phantom{\cdots}\cdots\phantom{\cdots}& +&b_{n-2}\, x_{n-2}\,\xi_1^2& +&\beta_1 \,\xi_1^3 &+&\beta_2\,\xi_1^2\,\xi_2&+\\
+& c_1 \,x_1\, \xi_2^2& +& \phantom{\cdots}\cdots\phantom{\cdots}& +&c_{n-2}\, x_{n-2}\,\xi_2^2& +&\gamma_1\, \xi_1\,\xi_2^2
&+&\gamma_2\,\xi_2^3&+\\ +& d_1 \,x_1\, \xi_1\,\xi_2& +& \phantom{\cdots}\cdots\phantom{\cdots}& +&d_{n-2}\, x_{n-2}\,\xi_1\,\xi_2& +&&&&\\
+&\rlap{\rm (terms of degree $3$ in $(x, \xi)$ and degree $\ge 2$ in $x$)} &&&&&&&& &\\
+&\rlap{$f_{40} \,\xi_1^4 \quad+\quad f_{31}\, \xi_1^3 \,\xi_2 \quad+\quad f_{22} \,\xi_1^2 \,\xi_2^2 
\quad+\quad f_{13} \,\xi_1 \,\xi_2^3 \quad+\quad f_{04} \,\xi_2^4\quad+\quad$} &&&&&&&& &\\ 
+&\rlap{\rm (terms of degree $4$ in $(x, \xi)$ and degree$\ge 1$ in $x$)} &&&&&&&& &\\ 
+&\rlap{\rm (terms of degree $\ge 5$ in $(x, \xi)$)}&&&&&&&& &\\ 
&&&&&&&&& &
\end{array}
$$
\caption{The coefficients of a formal power series in $(x, \xi)$}\label{table:coeff}
\end{table}
The following Proposition follows immediately from Lemma~\ref{lem:trans}.
\begin{proposition}\label{prop:AAeven}
Suppose that the rank of $H_{\phi, p}$ is $n-2$.
Then  the following condition is 
independent of the choice of $\phi$ and  $(x, \xi)$:
at least one of the coefficients
$\alpha_1, \alpha_2, \beta_1, \beta_2, \gamma_1, \gamma_2$ 
is not zero.
\end{proposition}
\begin{definition}\label{def:AAeven}
We say  that the singularity of $D$ at $p$ is \emph{of type $\AAeven$} if 
$H_{\phi, p}$ is of rank $n-2$ and 
at least one of the coefficients
$\alpha_1, \alpha_2, \beta_1, \beta_2, \gamma_1, \gamma_2$ 
is not zero.
\end{definition}
\begin{proposition_definition}\label{propdef:NCNR}
{\rm (1)}
The following three conditions are equivalent:
\begin{itemize}
\item[(i)] $\dim_k \OOO_{C(D, p), p}\le 4$,
\item[(ii)] $\OOO_{C(D, p), p}\cong k[[s, t]]/(s^2, t^2)$, and
\item[(iii)] the rank of the Hessian $H_{\phi, p}$ is $n-2$, and
\begin{equation}\label{eq:rank1}
\rank \left[
\begin{array}{ccc}
\alpha_1 & \beta_1 & \beta_2 \\
\alpha_2 & \gamma_1 & \gamma_2
\end{array}
\right]
=2.
\end{equation}
\end{itemize}
We say that the singularity of $D$ at $p$ is \emph{of type $\typeC$} if these conditions are satisfied.
\par
{\rm (2)}
We put
$$
T_1:=\sum_{i=1}^{n-2} d_i b_{\tau (i)} +f_{31}\quand 
T_2:=\sum_{i=1}^{n-2} d_i c_{\tau (i)} +f_{13}.
$$
The following three conditions are equivalent:
\begin{itemize}
\item[(i)] $\dim_k \OOO_{R(D, p), p}\le 4$,
\item[(ii)] $\OOO_{R(D, p), p}\cong k[[s, t]]/(s^2, t^2)$, and
\item[(iii)] the rank of the Hessian $H_{\phi, p}$ is $n-2$, and
\begin{equation}\label{eq:rank2}
\rank \left[
\begin{array}{cccc}
\alpha_1 & \beta_1 & \beta_2 & T_1 \\
\alpha_2 & \gamma_1 & \gamma_2 & T_2
\end{array}
\right]
=2.
\end{equation}
\end{itemize}
We say that the singularity of $D$ at $p$ is \emph{of type $\typeR$} if these conditions are  satisfied.
\end{proposition_definition}
\begin{remark}
As a corollary of Proposition-Definition~\ref{propdef:NCNR},
we see that the conditions~\eqref{eq:rank1} and~\eqref{eq:rank2} do not depend 
on the choice of the admissible formal parameter system $(x,\xi)$. 
This fact can also be proved directly by means of  Lemma~\ref{lem:trans}.
\end{remark}
For the proof of Proposition-Definition~\ref{propdef:NCNR},
we use the following easy lemma:
\begin{lemma}\label{lem:gs}
Let $g_1, \dots, g_l$ be degenerate formal power series   in  $k[[s, t]]$,
and let $J$ be the ideal of $k[[s, t]]$ generated by $g_1, \dots, g_l$.
We denote  by $a_i$ and $b_i$ the coefficients of $s^2$ and $t^2$ in $g_i$, respectively;
$$
g_i\;=\;a_i s^2+b_i t^2 +(\hbox{\rm terms of degree $\ge 3$}).
$$
Then the following equivalence holds:
$$
\dim_k k[[s, t]]/J\le 4
\;\;
\Longleftrightarrow
\;\;
J=(s^2, t^2)
\;\;
\Longleftrightarrow
\;\;
\rank
\left[
\begin{array}{ccc}
a_1 &\hskip -6pt \dots\hskip -6pt& a_l \\
b_1 &\hskip -6pt\dots\hskip -6pt& b_l
\end{array}
\right]=2.
$$
\end{lemma}
\begin{proof}[Proof of Proposition-Definition~\ref{propdef:NCNR}]
Since the proof of the assertion (1) is similar to and simpler than that of (2),
we prove only the assertion (2).
For simplicity,
we put
$$
D_i \phi:=\Der{\phi}{x_i},
$$
and let $I$ denote the ideal of $(\OOO_{X,P})\comp$ generated by
$\phi, D_1\phi, \dots, D_n\phi$ and $\rho$.
Suppose that $\rank H_{\phi, p}=n-2l$,
and let $(x_1, \dots, x_n)$ be an admissible formal parameter system with respect to $\phi$.
Then the linear parts of $\phi$ and $D_i\phi$ are given as follows;
$$
\phi^{[1]}=0,
\quad
(D_i\phi)^{[1]}=\begin{cases}
x_{\tau (i)} & \text{\rm if $i\le n-2l$} \\
0 & \text{\rm if $i> n-2l$}.
\end{cases}
$$
By Corollary~\ref{cor:sqrtdetH}, we also have
$$
\rho^{[1]}=
\begin{cases}
0 & \text{\rm if $n-2l\le n-4$} \\
d_1 x_1 +\dots+d_{n-2} x_{n-2} & \text{\rm if $ n-2l=n-2$}.
\end{cases}
$$
Therefore the classes of  $1\in k$ and $ x_{n-2l+1}, \dots, x_{n}$ 
span a $(2l+1)$-dimensional linear subspace in $(\OOO_{R(D,p),p})\comp=k[[x_1, \dots, x_n]]/I$.
Consequently the condition that $\dim_k \OOO_{R(D, p), p}\le 4$
implies  that $\rank H_{\phi, p}$ is $n-2$.
Hence it suffices to 
show the equivalence of the conditions (i), (ii) and (iii)
under the assumption that $\rank H_{\phi, p}=n-2$.
\par
\medskip
We use the admissible formal parameter system $(x, \xi)$  and Table~\ref{table:coeff}.
Let $I\sprime$ be the ideal of $k[[x, \xi]]$
generated by $D_1 \phi, \dots, D_{n-2} \phi$, and 
put
$$
s:=\xi_1\bmod I\sprime, \quad t:=\xi_2\bmod I\sprime.
$$
For $i=1, \dots, n-2$, we have
$$
D_i \phi\;=\; x_{\tau(i)} \,+\, b_i \xi_1^2 \,+\, c_i \xi_2 ^2 \,+\, d_i \xi_1 \xi_2 \;+\;
\left(\;\parbox{6cm}{\rm terms of degree $\ge 2$ in $(x,\xi)$ and degree $\ge 1$ in $x$, or degree $\ge 3$ in $(x,\xi)$}\;\right).
$$
Hence we obtain 
\begin{equation}\label{eq:xi}
x_i \bmod I\sprime = b_{\tau(i)} s^2 + c_{\tau(i)} t^2 + d_{\tau (i)} st + (\hbox{\rm terms of degree $\ge 3$ in $(s, t)$})
\end{equation}
for $i=1, \dots, n-2$.
In particular, we have $k[[x, \xi]]/I\sprime=k[[s, t]]$.
Using~\eqref{eq:xi},
we have
$$
\begin{array}{cclcc}
g_1  &:=& \phi \bmod I\sprime &=& \alpha_1 \, s^2 \,+\, \alpha_2\, t^2\,+\,(\hbox{\rm terms of degree $\ge 3$ in $(s, t)$}), \\
g_2  &:=& (D_{n-1} \phi) \bmod I\sprime &=& \beta_1 \, s^2 \,+\, \gamma_1\, t^2\,+\,(\hbox{\rm terms of degree $\ge 3$ in $(s, t)$}), \\
g_3  &:=& (D_{n} \phi)\bmod I\sprime &=& \beta_2 \, s^2 \,+\, \gamma_2\, t^2\,+\,(\hbox{\rm terms of degree $\ge 3$ in $(s, t)$}).
\end{array}
$$
On the other hand,
using~\eqref{eq:xi} again,
we see that
$$
\DDer{\phi}{\xi_1}{\xi_2}\,\bmod I\sprime, 
\quad
\DDer{\phi}{x_i}{\xi_1}\,\bmod I\sprime,
\quand
\DDer{\phi}{x_i}{\xi_2}\,\bmod I\sprime
$$ 
are equal to 
$$
\begin{array}{l}
 f_{31} \,s^2 + f_{13}\,t^2 +\dsum_{i=1}^{n-2} d_i \,(b_{\tau(i)} s^2 + c_{\tau(i)} t^2 + d_{\tau (i)} st)
+(\hbox{\rm terms of degree $\ge 3$ in $(s, t)$}),\\
d_i \,t +(\hbox{\rm terms of degree $\ge 2$ in $(s, t)$}) \mystrutd{7pt}\quand\\
d_i \,s +(\hbox{\rm terms of degree $\ge 2$ in $(s, t)$}),
\end{array}
$$
respectively.
Moreover, we have
$$
\DDer{\phi}{x_i}{x_j}\,\bmod I\sprime\;=\;
\begin{cases}
1+\,(\hbox{\rm terms of degree $\ge 1$ in $(s, t)$}) & \text{if $i=\tau (j)$,} \\
\phantom{1+\,\,\,}\,(\hbox{\rm terms of degree $\ge 1$ in $(s, t)$})& \text{if $i\ne \tau (j)$,}
\end{cases}
$$
for $i, j=1, \dots, n-2$.
By Corollary~\ref{cor:sqrtdetH} and using  $\sum_{i=1}^{n-2} d_i d_{\tau(i)}=0$, we obtain 
$$
g_4\;:=\;\rho\bmod I\sprime \;=\;T_1 \, s^2\, +\,
  T_2 \, t^2\,+\,(\hbox{\rm terms of degree $\ge 3$ in $(s, t)$}).
$$
Since
$$
(\OOO_{R(D, p), p})\comp\;\;=\;\;k[[x, \xi]]/I\;\; =\;\; k[[s, t]]/( g_1, g_2, g_3, g_4),
$$
we obtain the equivalence of the conditions (i)-(iii) from Lemma~\ref{lem:gs}. 
\end{proof}
The following remark is crucial in the proof of Proposition~\ref{prop:eval}.
\begin{remark}\label{rem:open}
Each of the conditions that 
the singularity of $D$ at $p$ is of type $\AAeven$, $\typeC$ and  $\typeR$ is 
an \emph{open} condition on the coefficients of the formal power series $\phi$
with  $\rank H_{\phi, p}\le n-2$ in the following sense.
We denote by $\mmm$ the maximal ideal of $k[[x_1, \dots, x_n]]$.
For a positive integer $d\ge 2$, we put
$$
V_d:=\mmm^2/\mmm^{d+1}.
$$
For $\phi\in \mmm^2$,
let $\bar \phi\in V_d$ denote the class of $\phi$ modulo $\mmm^{d+1}$.
By Lemma~\ref{lem:irredfm},
the locus
$$
\Delta:=\set{\bar\phi\in V_2}{\rank H_{\phi,p}\le n-2}
$$
is a closed  irreducible hypersurface of $V_2$.
For $d\ge 3$, 
let $\Delta _d$ denote  the pull-back of
$\Delta$ by the natural linear homomorphism
$V_d\to V_2$.
Note that the types
$\AAeven$, $\typeC$ and $\typeR$ of the hypersurface singularities are 
defined by using   the coefficients of terms of degree $\le 3$, $\le 3$ and $\le 4$
of the formal power series defining the hypersurface, respectively.
Hence
the loci
$$
\renewcommand{\arraystretch}{1.28}
\begin{array}{rcl}
\Delta \AAeven &:=& \set{\bar\phi\in V_3}{\hbox{the singularity of $\phi=0$ at $p$ is of type $\AAeven$}}, \\
\Delta \typeC &:=& \set{\bar\phi\in V_3}{\hbox{the singularity of $\phi=0$ at $p$ is of type $\typeC$}},\\
\Delta \typeR &:=& \set{\bar\phi\in V_4}{\hbox{the singularity of $\phi=0$ at $p$ is of type $\typeR$}},
\end{array}
$$
are well defined.
The loci $\Delta \AAeven $ and $\Delta \typeC$ are Zariski open dense subsets of $\Delta_3$,
and the locus $\Delta \typeR$ is a Zariski open dense subset of $\Delta_4$.
\end{remark}
\section{The  singular scheme  $\Gamma_{m-1}$}\label{sec:Gamma}
We define a  homomorphism
$$
\map{\gamma^*}{k [[\sixvars]]}{ k[[\threevars]]}
$$
over $k$ by
$$
\gamma ^* u_1 =x_1^2,
\quad
\gamma ^* v_1 =x_1 z,
\quad
\gamma ^* u_2 =x_2^2,
\quad
\gamma ^* v_2 =x_2 z,
\quad
\gamma ^* w =x_1 x_2 z,
\quad
\gamma ^* y =z,
$$
and denote by
$$
\Gamma:=\Spec (k[[\sixvars]]\,/\Ker \gamma^* )
$$
the scheme-theoretic image of the morphism
$$
\map{\gamma}{\Spec k[[\threevars]]}{\Spec k [[\sixvars]]}
$$
induced from $\gamma^*$.
The critical subscheme of $\gamma$ is defined  by $z^2=0$.
Calculating the Gr\"obner basis (see~\cite{MR1417938}) of the ideal 
$$
\langle \;
u_1 -x_1^2, \;
v_1 -x_1 z, \;
u_2 -x_2^2, \;
v_2 -x_2 z, \;
w -x_1 x_2 z, \;
y -z \;
\rangle
$$
in the polynomial ring 
$k[\threevars, \sixvars]$ under  the pure lexicographic order
$$
x_1 \;>\; x_2 \;>\; z \;>\; w \;>\; v_2 \;>\; v_1 \;>\; y \;>\; u_2 \;>\; u_1, 
$$
we see that 
$\Gamma$ is  a $3$-dimensional singular scheme  defined by the equations~\eqref{eq:defeqGamma} 
given in Introduction.
The singular locus of $\Gamma$ is equal to the $2$-dimensional plane
$$
v_1=v_2=w=y=0.
$$
For a $k$-algebra $R$, we denote by
$$
\map{\gamma^R}{\Spec R[[\threevars]]}{\Spec R [[ \sixvars]]}
$$
the morphism obtained from $\gamma$ by $k\inj R$,
and by
$$
\Gamma^R:=\Spec (R[[\sixvars]]\,/\,\Ker \,(\gamma^R)^* ),
$$
the scheme-theoretic image of $\gamma^R$.
\begin{theorem}\label{thm:Gamma}
Let $V$ be a smooth variety of dimension $\dimV \ge 3$,
and  let 
$\shortmap{\phi}{V}{W}$ 
be a morphism to a smooth variety $W$.
Suppose that 
there exists a smooth hypersurface $H$ of $V$ such that,
at every point $P$ of $H$, the kernel of
$$
\map{d_P \phi}{T_P (V)}{T_{\phi(P)} (W)}
$$
is of dimension $2$ and is contained in $T_P (H)$.

Let $P$ be a point of $H$.
We denote by $R$ the ring of formal power series
$k[[ x_3, \dots, x_{\dimV-1}]]$ of $\dimV-3$ variables.
Then there exists an isomorphism
$$
\mapisom{\iota}{(V, P)\comp}{\Spec R[[\threevars]]}
$$
over $k$ with the following properties:
\begin{itemize}
\item[(i)]
the divisor $(H, P)\comp$ of $(V, P)\comp$ is mapped by the isomorphism $\iota$ to the divisor given by $z=0$,
and
\item[(ii)]
there exists a unique morphism $\psi : \Gamma^R\to (W, \phi(P))\comp $
such that 
the formal completion $\phi_P\comp$ of $\phi$ at $P$ factors  as 
the following diagram:
\begin{equation*}\label{eq:diag}
\begin{array}{ccc}
(V, P)\comp & \maprightsp{\phi_P\comp}& (W, \phi(P))\comp \\
\mapdownleftright{\iota}{\wr}&&\mapupright{\psi} \\
\Spec R[[\threevars]] & \maprightsb{\gamma^R}& \phantom{a}\Gamma^R.
\end{array}
\end{equation*}
\end{itemize}
\end{theorem}
\begin{proof}
We denote by
$$
\map{\phi_H}{H}{W}
$$
the restriction of $\phi$ to $H$.
The kernel $\KKK$ of the homomorphism
$$
\map{d \phi_H}{\Tan{H}}{\phi_H^* \Tan{W}}
$$
is an integrable  subbundle of $\Tan{H}$ with rank $2$
by the assumption.
Let $P$ be an arbitrary  point of $H$.
By~\cite[Proposition 6]{MR0450263Exp6}, 
there exist a basis $(D_1, D_2)$ of the $\OOO_{H, P}$-module  $\KKK\otimes \OOO_{H, P}$ 
and a system of uniformizing variables $(y_1, \dots y_{\dimV-1})$ of $H$ at $P$
such that
$$
D_i (y_j)=\delta_{ij}\qquad(i=1, 2, \;\;j=1, \dots, \dimV-1).
$$
In other words, the submodule $\KKK\otimes \OOO_{H, P}$  of $\Tan{H} \otimes \OOO_{H, P}$
is equal to
$$
\OOO_{H, P} \Der{}{y_1}\, \oplus \, \OOO_{H, P} \Der{}{y_2}.
$$
We can choose a formal parameter system $(x_1, \dots, x_{\dimV-1}, z)$ of $V$ at $P$
such that $H$ is defined by $z=0$,
and that the restriction $x_i\rest H$ of $x_i$ to $H$ is equal to $y_i$.
Then, for every element $g\in (\OOO_{W,\phi(P)})\comp$,
we have
$$
\Der{\;(\phi_H^* g)}{y_1} \;=\;\Der{\;(\phi^* g)}{x_1}\,\left|\,\atop{}{H}\right. \;\equiv\;  0
\quand
\Der{\;(\phi_H^* g)}{y_2} \;=\;\Der{\;(\phi^* g)}{x_2}\,\left|\,\atop{}{H}\right. \;\equiv\;  0;
$$
that is, $\phi^* g $ is contained in  the  subring $A$ of $ (\OOO_{V,P})\comp=k[[x_1, \dots, x_{\dimV-1}, z]]$ defined by
$$
A:=\left\{ f\in (\OOO_{V,P})\comp\;\;\left|\;\;\hbox{$\dDer{f}{x_1}$ and $\dDer{f}{x_2}$ are contained in 
the ideal $(z)$}\right.\right\}.
$$
Putting
$R:=k[[ x_3, \dots, x_{\dimV-1}]]$,
we have
$$
A\;\;= \;\;R[[x_1^2, x_2^2, z]] \;+\; z\,  R[[\threevars]]\;\;\; \st \;\;\; R[[\threevars]]\;=\; (\OOO_{V,P})\comp.
$$
Therefore, the subring $A$ is the image of the homomorphism
$$
\map{(\gamma^R )^*}{R[[\sixvars]]}{R[[\threevars]]}.
$$
We denote by
$$
\mapisom{j}{R[[\sixvars]]\,/\,(\Ker\,(\gamma^R)^*)}{A}
$$
the induced isomorphism.
Let $(t_1, \dots, t_{q})$ be a formal parameter system of $W$ at $\phi (P)$, where $q=\dim W$.
Since $\phi^* t_i\in A$,
there exists a unique element 
$$
\psi _i \;\in \;R[[\sixvars]]\,/\,(\Ker\,(\gamma^R)^*)
$$
such that $j (\psi_i)=\phi^* t_i$.
Define a morphism  
$$
\map{\psi}{\Gamma^R}{(W, \phi(P))\comp}
$$
by $\psi^* t_i=\psi_i$.
Then $\phi_P\comp$ factors as $\psi\circ \gamma^R\circ\iota$,
where $\iota$ is the isomorphism given by the canonical isomorphism 
$R[[\threevars]]\isom  (\OOO_{V, P})\comp$.
\end{proof}
\begin{definition}
When $R$ is the ring of  formal power series in $(m-4)$ variables with coefficients in $k$,
the scheme-theoretic image  $\Gamma^R$ is an $(m-1)$-dimensional singular scheme,
which we will denote by $\Gamma_{m-1}$.
\end{definition}
\begin{remark}
A  normal form theorem similar to Theorem~\ref{thm:Gamma}
for a morphism from a smooth surface to a smooth curve was proved in~\cite{MR1157319}.
\end{remark}
\section{The discriminant variety}\label{sec:reminder}
In this section, we recall the definitions and results in \cite[Section 3]{char3},
which are valid in any characteristics.
\par
\medskip
Let $\XXX$ be a projective variety of dimension $n>0$,
$\LLL$ a line bundle on $\XXX$,
and $M$ a linear subspace of $H^0 (\XXX, \LLL)$ with dimension $m+1\ge 2$.
We denote by
$$
\PM:=\P_* (M)
$$
the parameter space of 
the $m$-dimensional linear system $|M|$ of divisors on $\XXX$ corresponding to $M$,
and put
$$
\XX:=\XXX\sm (\Sing (\XXX) \cup \Bs (|M|)),
$$
where $\Sing (\XXX)$ is the singular locus of $\XXX$ and $\Bs (|M|)$ is the base locus of $|M|$.
We denote by $\Psi: \XX\to\PM\dual$
the morphism induced by $\PM$, and set
$$
\X:=\set{p\in \XX}{\hbox{the homomorphism $d_p \Psi: T_p (\XX)\to T_{\Psi (p)} (\PM\dual) $ is injective}}.
$$
Note that $\X=\XXX$ if $\XXX$ is smooth and $|M|$ is very ample.
\begin{assumption}
Throughout the paper,
we assume that $m>n$,  and that $\X$ is dense in $\XXX$.
\end{assumption}
\par
\medskip
For a non-zero vector $f\in M$,  we denote by $[f]\in \PM$ the corresponding point of $\PM$, 
by $\DDDf$ the divisor of $\XXX$ defined by $f=0$, and  by $\DDf$ the intersection $\DDDf\cap\XX$.
Let 
$$
\pr_1 : \XXP \to X\quand \pr_2 :\XXP \to\PM
$$
be the projections.
We put
$$
\widetilde {\LLL}:=\pr_1\sp * \LLL \otimes \,\pr_2\sp * \OOO_{\PM} (1).
$$
There exists a canonical isomorphism
$M\dual \isom H^0 (\PM, \OOO_{\PM}(1))$
unique up to multiplicative constants.
Combining this isomorphism with the inclusion $M\inj H^0(X, \LLL)$,
we obtain a natural  homomorphism
$$
\Hom (M, M) =M\otimes M\dual \;\; \to\;\; H^0 (\XXP, \widetilde {\LLL}).
$$
We denote by 
$$
\sigma\in H^0 (\XXP, \widetilde {\LLL})
$$ 
the global section of $\widetilde {\LLL}$ corresponding to the identity homomorphism of $M$.
Note that $\sigma$ is determined uniquely up to multiplicative constants.
The zero locus of $\sigma$ is equal to
$$
\set{\pf\in \XXP}{p\in \DDf}.
$$
Let $\fDD$ be the subscheme of $\XXP$ defined by $\sigma=0$, and 
let 
$$
p_1 :\fDD \to \XX\quand p_2 :\fDD \to \PM
$$
be the projections.
Then $\fDD$ is smooth of dimension $m+n-1$,
and 
$p_2 :\fDD\to\PM$ is the universal family of the divisors $\DDf\;([f]\in \PM)$.
We have a natural homomorphism
$$
\map{d\sigma}%
{\Tan{\XXP}\otimes\OOO\sb{\fDD}}%
{\wt{\LLL}\otimes\OOO\sb{\fDD}}
$$
defined by $\Theta\mapsto (\Theta\sigma )\rest \fD$,
where $\Theta$ is a local section of $\Tan{\XXP}$ considered as a derivation of $\OOO\sb{\XXP}$.
We denote by
$$
\map{d\sigma\sb{\XX}}{p_1^* \Tan{\XX}}%
{\wt{\LLL}\otimes\OOO\sb{\fDD}}
$$
the restriction of $d\sigma$ to the direct factor $p_1^* \Tan{\XX}$
of 
$$
\Tan{\XX\times\PM}\otimes\OOO\sb{\fDD}=p_1^* \Tan{\XX}\oplus p_2^* \Tan{\PM}.
$$
Let $\fCC$ be the critical subscheme of $p_2: \fDD\to\PM$.
\begin{proposition}\label{prop:fCC}
{\rm (1)}
The  scheme $\fCC$ is equal to the degeneracy subscheme  of $d\sigma\sb{\XX}$.
\par
{\rm (2)}
The intersection $\fC$ of $\fCC$ and $\XP$ is smooth, irreducible and  of dimension $m-1$.
\end{proposition}
\begin{corollary}
{\rm (1)}
The support of $\fCC$ is equal to the set 
$$
\set{\pf\in \fDD}{\hbox{\rm $\DDf$ is singular at $p$}}.
$$
\par
{\rm (2)}
Let $P=\pf$ be a point of $\fCC$,
and let $\pi_2\sprime :\fCC\to\PM$ be the projection.
Then,
as  a subscheme of $(\pr_2\inv ([f]), P)\comp\cong (X, p)\comp$, 
the completion $(\pi_2\sp{\prime-1} ([f]), P)\comp $ of the fiber of $\pi_2\sprime$
at $P$ coincides with $C(D_{[f]}, p)$ defined in Definition~\ref{def:CR}.
\end{corollary}
We denote by  
$$
\pi_1 :\fC \to \X\quand \pi_2 :\fC \to \PM
$$
the projections from $\fC=\fCC\cap (\XP)$.
We can define 
the \emph{universal Hessian}
$$
\map{\fH}{ \pi_1^* \Tan{\X} \otimes\pi_1^* \Tan{\X} }{ \wt{\LLL}\otimes \OOO\sb{\fC}}
$$
on $\fC$. 
Because we are in characteristic $2$,
the universal Hessian $\fH$ on $\fC$  is \emph{anti}-symmetric.
Let
$$
\map{\fHt }{\pi_1^* \,\Tan{\X} }{ \wt{\LLL}\otimes (\pi_1^* \,\Tan{\X})\dual }
$$
be the homomorphism obtained from $\fH$.
Let $\fE$ be the critical subscheme of the projection $\pi_2 : \fC \to \PM$.
\begin{proposition}\label{prop:fE}
The  scheme $\fE$ 
is equal to the degeneracy subscheme   of $\fHt$.
\end{proposition}
\begin{corollary}
The support of $\fE$ is equal to the set of all point $(p, [f])\in \fC$
such that the Hessian of the hypersurface singularity $p\in D_{[f]}$ is degenerate.
\end{corollary}
Since the rank of an anti-symmetric form is always even, 
we obtain the following:
\begin{corollary}\label{prop:odd}
Suppose that $n$ is odd.
Then $\fE$ coincides with $\fC$.
\end{corollary}
\section{Main results}\label{sec:evenmain}
From now on until the end of the paper, we assume that  $n$ is even.
\begin{definition}\label{def:RRR}
Let $\fR\subset \fC$ be the degeneracy subscheme $\Degasymm (\fH)$ of the anti-symmetric form  $\fH$,
and let
$$
\varpi_1: \fR\to \X\quand \varpi_2: \fR\to \PM
$$
be the projections.
\end{definition}
\begin{proposition}\label{prop:irred}
{\rm (1)}
The subscheme $\fR$ is of codimension $\le 1$ in $\fC$.
If $\fR$ is of codimension $1$ in $\fC$, then 
$\fE$ coincides with $2\,\fR$ as a divisor of $\fC$.
\par
{\rm (2)}
Let $P=\pf$ be a point of $\fR$.
Then,
as  a subscheme of $(\pr_2\inv ([f]), P)\comp\cong (X, p)\comp$, 
the completion $(\varpi_2\inv ([f]), P)\comp $ of the fiber of $\varpi_2$
at $P$ coincides with $R(D_{[f]}, p)$ defined in Definition~\ref{def:CR}.
\end{proposition}
\begin{proof}
The  assertion (1) follows from  Remark~\ref{rem:as_a_divisor} and  Proposition~\ref{prop:fE}.
The  assertion (2) follows from the definition.
\end{proof}
It is possible that 
 $\fR$ is of codimension $0$ in $\fC$.
The following is a classical example due to Wallace~\cite{MR0080354}.
See also~\cite{MR1794260}.
\begin{example}\label{example:Fermat1}
Let $X$ be the Fermat hypersurface of degree $2^\nu+1$ in $\P^{n+1}$
with $\nu\ge 1$,
and let $|M|$ be the complete linear system $|\OOO_{X} (1)|$.
Then the rank of the Hessian of $D_{[f]}$ at $p$  is zero at every point $\pf$ of $\fC=\fCC$.
In particular, $\fR$ coincides with $\fC$.
\end{example}
Let $P=\pf$ be a point of $\fR$,
so that $D_{[f]}$ is singular at $p$ and the Hessian of $D_{[f]}$ at $p$ is of rank $\le n-2$.
We investigate the formal completion of the morphism $\pi_2 : \fC \to \PM$ at $P$.
For this purpose, 
we introduce a good formal parameter system of $\XP$ at $P$.
\par
\medskip
Because $p\in \X$, we can choose a basis $\vb_0, \dots, \vb_m$ of the vector space
$M$ in such a way that the following conditions (i)-(iv)
hold. Let $\tvb_i$ be the global section of $\LLL$ corresponding to $\vb_i$,
and let $D_i$ be the divisor defined by $\tvb_i=0$.
Then
\begin{itemize}
\item[(i)] $\vb_0=f$, so that $D_0=\DDf$ passes through $p$ and is singular at $p$, 
\item[(ii)] $p\in D_i$ for $i=0, 1, \dots, m-1$, and $p\notin D_m$, 
\item[(iii)] $D_1, \dots, D_n$ are smooth at $p$ and intersect transversely at $p$, and 
\item[(iv)] $D_{n+1}, \dots, D_{m-1}$ are singular at $p$.
\end{itemize}
For $i=0, \dots, m-1$, we put
$$
\phi_i := \tvb_i/\tvb_m,
$$
which is a rational function on $X$ regular at $p$.
Then $(\phi_1, \dots, \phi_n)$ forms a local  parameter system of $X$ at $p$ by (iii).
We put
$$
x_i :=\phi_i\quad (i=1, \dots, n).
$$
We will regard $(x_1, \dots, x_n)$ as a formal parameter system of $X$ at $p$, and consider 
$\phi_0$ and $\phi_{n+1}, \dots, \phi_{m-1}$ as formal power series of $(x_1, \dots, x_n)$,
which have no terms of degree $\le 1$ by (i) and (iv).
Note that  the rank of the anti-symmetric  matrix $H_{\phi_0, p}$ is $\le n-2$,
because $P\in \fR$.
Using Remark~\ref{rem:make_admissible},
we can further assume the following by a linear transformation of $\vb_1, \dots, \vb_n$:
\begin{itemize}
\item[(v)] the formal parameter system $(x_1, \dots, x_n)$ is admissible with respect to $\phi_0$.
\end{itemize}
Let $(Y_0, \dots, Y_m)$ be the linear coordinates of $M$
with respect to the basis $\vb_0, \dots, \vb_m$,
and let
$(y_1, \dots, y_m)$ be the affine coordinate system of $\PM$ given by
$$
y_i:=Y_i/Y_0.
$$
Then the point $[f]\in \PM$ is the origin.
We regard 
$$
(x_1, \dots, x_n, y_1, \dots, y_m)
$$
as a formal parameter system of $\XP$ at $P=\pf$.
\par
\medskip
We can regard the coordinates $Y_i$  as global sections of $\OOO_{\PM} (1)$.
Recall that $\sigma$ is the canonical section of
$\wt{\LLL}:=\pr_1\sp * \LLL\otimes \pr_2\sp * \OOO_{\PM} (1)$.
Multiplying $\sigma$ by a suitable non-zero constant,
we have 
\begin{equation*}\label{eq:sigma_explicit}
\sigma=\sum\sb{i=0}\sp {m} \tvb_i \otimes Y_i.
\end{equation*}
We put
$$
\Phi:=\left(\sum\sb{i=0}\sp {m} \tvb_i \otimes Y_i \right)/(\tvb_m \otimes Y_0),
$$
which is a rational function of $\XP$ regular at $P$.
By definition, we have
$$
\Phi=\phi_0 +y_1x_1 + \cdots +y_n x_n + y_{n+1} \phi_{n+1} + \cdots + y_{m-1}\phi_{m-1} +y_m.
$$
The scheme  $\fDD$ is defined by $\Phi=0$ locally at $P$.
For simplicity,
we put
$$
D_i \Phi:= \Der{\Phi}{x_i}.
$$
By Proposition~\ref{prop:fCC},
the subscheme $\fC$ of $\XP$  is defined by
$$
\Phi=D_1 \Phi=\dots=D_n \Phi=0
$$
locally at $P$.
We put
$$
R:=\;\;f_n \left(\dots,\;\;\DDer{\Phi}{x_i}{x_j},\;\; \dots
\right)\;\;=\;\;\sqrt{\det\left(\DDer{\Phi}{x_i}{x_j}\right)}\;\;,
$$
where $f_n$ is the polynomial defined in~\S\ref{subsec:anti-symmetric}.
The subscheme  $\RRR$  of $\XP$ is defined by
$$
\Phi=D_1 \Phi=\dots=D_n \Phi=R=0
$$
locally at $P$.
\par
\medskip
Note that the divisor $\DDf$ of $X$ is defined by $\phi_0=0$ locally at $p$.
Let $2r$ be the rank of the Hessian  $H_{\phi_0, p}$ of $\DDf$ at $p$.
Since $(x_1, \dots, x_n)$ is admissible with respect to $\phi_0$,
we see that the linear parts of the formal power series $\Phi$ and $D_i\Phi$ are as follows:
\begin{equation}\label{eq:Philin}
\renewcommand{\arraystretch}{1.2}
\begin{array}{lcll}
\Phi^{[1]} &=& y_m, &   \\
(D_i \Phi)^{[1]} &=& x_{\tau(i)} +y_i & \textrm{for $i=1, \dots, 2r$,}  \\
(D_i \Phi)^{[1]} &=& y_i & \textrm{for $i=2r+1, \dots, n$.} 
\end{array}
\end{equation}
When $2r=n-2$,
we use the  formal parameter system $(x_1, \dots, x_{n-2}, \xi_1, \xi_2)$
admissible with respect to $\phi_0$
as in \S\ref{sec:HS},
and adopt  Table~\ref{table:coeff} in order to denote the coefficients of 
$\phi_0$.
We also denote by $e_j\in k$ the coefficients
of the term $\xi_1\xi_2$ in $\phi_j$ for $j=n+1, \dots, m-1$.
By Corollary~\ref{cor:sqrtdetH},  we see that the linear part of $R$ is equal to 
\begin{equation}\label{eq:Rlin}
\begin{array}{rcl}
R ^{[1]} &=& 
\begin{cases}
0 & \textrm{if $2r<n-2$, }\\
\left(\dDDer{\Phi}{\xi_1}{\xi_2} \right)\sp{[1]} & \textrm{if $2r=n-2$, }\\
\end{cases}
\mystrutd{25pt}
\\
&=& 
\begin{cases}
0 & \textrm{if $2r<n-2$, }\\
\displaystyle {d_1 x_1 +\dots+ d_{n-2}x_{n-2} +\phantom{aaaa} }\mystruth{12pt}
\atop \displaystyle { \phantom{a}+ e_{n+1} y_{n+1} + \dots+ e_{m-1}y_{m-1} \mystruth{9pt}}  & \textrm{if $2r=n-2$. }\\
\end{cases}
\\
\end{array}
\end{equation}
\subsection{A normal form theorem}
We put
$$
\fRsm:=\set{P\in \fR}{\hbox{$\fR$ is smooth of dimension $m-2$ at $P$}},
$$
and denote by
$\shortmap{\vpism_2}{\fRsm}{\PM}$
the projection.
\par
\medskip
\begin{theorem}\label{thm:evenmain}
Let $P=\pf$ be a point of $\fR$.
\par
{\rm (1)}
If $P\in \fRsm$, then the rank of the Hessian $H_{\phi_0, p}$ of $\DDf$ at $p$ is $n-2$.
Conversely,
suppose that $H_{\phi_0, p}$ is of rank $n-2$.
Then $P$ is a point of $\fRsm$ if and only if
at least one of 
$$
d_1, \dots, d_{n-2}, e_{n+1}, \dots, e_{m-1}
$$
is not zero.
\par
{\rm (2)}
The kernel $\KKK$ of the homomorphism
$$
\map{d\vpism_2}{\Tan{\fRsm}}{(\vpism_2)^* \Tan{\PM}}
$$
is a subbundle of $\Tan{\fRsm}$ with rank $2$.
\end{theorem}
\begin{proof}
The assertion (1) follows immediately from~\eqref{eq:Philin} and~\eqref{eq:Rlin}.
Suppose that $P\in \fRsm$.
Then the kernel of $d_P \vpism_2: T_P (\fRsm)\to T_{[f]} (\PM)$ is of dimension $2$
and generated by the vectors
$$
\left(\Der{}{x_{n-1}}\right)_P\quand \left(\Der{}{x_n}\right)_P
$$
of $T_P (\XP)$.
Since this fact holds at every point $P$ of $\fRsm$,
we see that $\KKK$ is a subbundle of $\Tan{\fRsm}$ with rank $2$.
\end{proof}
\begin{corollary}\label{cor:empty}
Suppose that $(n, m)=(2, 3)$. Then $\fRsm$ is empty.
\end{corollary}
\par
\medskip
By definition, the kernel $\KKK$ of $d\vpism_2$ is an integrable subbundle of $\Tan{\fRsm}$.
From Proposition~\ref{rem:facquot}, we obtain the following:
\begin{corollary}\label{cor:factorKKK}
The morphism
$\vpism_2: \fRsm\to \PM$ factors through the quotient morphism
$\shortmap{q}{\fRsm}{(\fRsm)\sp{\KKK}}$
by  $\KKK$, which is finite of degree  $4$.
\end{corollary}
Combining Theorems~\ref{thm:Gamma} and~\ref{thm:evenmain},
we obtain the following  normal form theorem for the morphism $\pi_2:\fC\to\PM$ at a point $P=\pf$ of $\fRsm$.
\begin{corollary}\label{cor:ci}
Let $P=\pf$ be a point of $\fRsm$.
Then there exist an isomorphism
$$
\mapisom{\iota}{(\fC, P)\comp}{\Spec k[[w_1, \dots, w_{m-1}]]}
$$
and a  morphism
$$
\map{\psi}{\Gamma_{m-1}}{(\PM, [f])\comp}
$$
such that the formal completion
$$
\map{(\pi_2)_P\comp}{(\fC, P)\comp}{(\PM, [f])\comp}
$$
of $\pi_2$ at $P$ factors as $\psi\circ\gamma^R\circ\iota$,
where
$$
\map{\gamma^R}{\Spec k[[w_1, \dots, w_{m-1}]]}{\Gamma_{m-1}}
$$
is the morphism defined in~\S\ref{sec:Gamma}.
\end{corollary}
\subsection{General plane sections}
Let $P=\pf$ be a point of $\fR$.
We choose 
linear subspaces
$L$ and $\Lambda$ of $\PM$ with dimension $3$ and $2$, respectively, 
such that
$$
[f]\;\in\;\Lambda\;\st\;L.
$$
We set
$$
\SL:=\pi_2\inv (L)\;\st\;\fC\quand
\CL:=\pi_2\inv (\Lambda)\;\st\;\fC,
$$
and denote by
$$
\map{\pi_L}{\SL}{L}\quand
\map{\pi_\Lambda}{\CL}{\Lambda}
$$
the restrictions of $\pi_2:\fC\to\PM$ to $\SL$ and to $\CL$, respectively.
\begin{proposition}\label{prop:SL}
Suppose that $L$ is chosen generically over $k$.
The scheme  $\SL$ is smooth of dimension $2$ at $P$
if and only if the rank  of the Hessian $H_{\phi_0, p}$ of $D_{[f]}$ at $p$ is $n-2$.
\end{proposition}
\begin{proof}
We can assume that $L$ is defined by 
\begin{equation}\label{eq:abc}
y_i\;=\;A_i\, y_{n-1} \,+\, B_i \,y_n \,+\,C_i \,y_m\qquad (i\ne n-1, n, m),
\end{equation}
where $A_i, B_i,C_i$ are generic over  $k$.
Then the assertion follows from the linear parts~\eqref{eq:Philin}
of the defining equations of $\fC$.
\end{proof}
\begin{proposition}\label{prop:cut}
Suppose that $L$ and $\Lambda$ are chosen generically over $k$,
and that $\SL$ is smooth of dimension $2$ at $P$.
\par
{\rm (1)}
If $\CL$ is of codimension $1$ in $\SL$,
then $\CL$ has a degenerate singularity at $P$.
\par
{\rm (2)}
The following two conditions are equivalent:
\begin{itemize}
\item[(i)]
$\CL$ is of codimension $1$ in $\SL$, and the multiplicity of $\CL$ at $P$ is $2$, 
\item[(ii)]
the singularity of  $\DDf$ at $p$ is  of type $\AAeven$.
\end{itemize}
\par
{\rm (3)}
Suppose  that $\CL$ is of codimension $1$ in $\SL$, that the 
multiplicity of $\CL$ at $P$ is $2$, 
and that $P$ is a point of $\fRsm$.
Then $P$ is an ordinary cusp of $\CL$ 
if and only if the singularity of  $\DDf$ at $p$ is   of type $\typeR$.
\end{proposition}
\begin{proof}
By Proposition~\ref{prop:SL} and the assumption, the rank of $H_{\phi_0, p}$ is $n-2$.
We use the formal parameter system  $(x_1, \dots, x_{n-2}, \xi_1, \xi_2)$
admissible with respect to $\phi_0$ as in \S\ref{sec:HS},
and refer to
Table~\ref{table:coeff} for the name of  coefficients of $\phi_0$.
We also use the defining equations~\eqref{eq:abc} of the linear subspace $L$.
By~\eqref{eq:Philin} and~\eqref{eq:abc}, we see that
$$
s:=\xi_1\rest\SL\quand t:=\xi_2\rest\SL
$$
form a formal parameter system of $\SL$ at $P$.
We put
$$
u:=y_{n-1}\rest {L}, \quad v:=y_n\rest {L}\quand w:=y_m\rest {L},
$$
which form an affine coordinate system of $L$ with the origin $[f]$.
For a formal power series $F$ of $(x_1, \dots, x_{n-2}, \xi_1, \xi_2,  y_1, \dots, y_m)$,
we denote by $F_L$ the formal power series of 
$$
(x_1, \dots, x_{n-2}, s, t, u, v, w)
$$
obtained from $F$ by making the substitutions
\begin{eqnarray*}
&& y_{n-1}=u, \quad y_n=v, \quad y_m=w, \\
&& y_i=A_i u + B_i v +C_i w\quad (i\ne n-1, n, m), \quand\\
&& \xi_1=s, \quad \xi_2=t.
\end{eqnarray*}
Regarding 
\begin{equation*}\label{eq:PhiL}
\Phi_L=(D_1 \Phi)_{L}=\cdots=(D_n\Phi)_{L}=0
\end{equation*}
as equations with indeterminates 
$x_1$, \dots, $x_{n-2}$,  $u$, $v$, $w$ and with coefficients in $(\OOO_{\SL, P})\comp=k[[s, t]]$, 
and solving them, 
we obtain the formal power series expansion 
$$
\renewcommand{\arraystretch}{1.28}
\begin{array}{ccccl}
x_i\rest\SL &=& X_i (s, t) & = & \sum X_{i, \mu\nu}\, s^\mu \,t^\nu \quad (i=1, \dots, n-2), \\
u\rest\SL &=& U (s, t) & = & \sum U_{\mu\nu} \,s^\mu \,t^\nu,\\
v\rest\SL &=& V (s, t) & = & \sum V_{\mu\nu} \,s^\mu \,t^\nu, \\
w\rest\SL &=& W (s, t) & = & \sum W_{\mu\nu} \,s^\mu \,t^\nu, \\
\end{array}
$$
of the functions 
$x_1|\SL$, \dots, $x_{n-2}|\SL$,  $u|\SL$, $v|\SL$, $w|\SL$ on $\SL$ 
in $(s, t)$.
The formal power series $U$, $V$ and $W$ give the formal completion of $\pi_L: \SL\to L$ at $P$.
We will calculate the homogeneous parts $U^{[d]}$, $V^{[d]}$, $W^{[d]}$ 
of degree $d$
of them  up to $d=3$.
From $\Phi_L=0$ and $(D_{n-1} \Phi)_L=(D_{n} \Phi)_L=0$,
we obtain
$$
U^{[1]}=V^{[1]}=W^{[1]}=0.
$$
From $(D_i\Phi)_L=0$ for $i=1, \dots, n-2$, we obtain
$$
X_i^{[1]}=0\qquad(i=1, \dots, n-2).
$$
Looking at the homogeneous parts of degree $2$ in $(s, t)$
of $\Phi_L$ and $(D_i \Phi)_{L}$,
we obtain the following equations:
$$
\renewcommand{\arraystretch}{1.2}
\begin{array}{lcll}
\alpha_1\, s^2 \,+\, \alpha_2\, t^2 \,+\, W^{[2]} &=&0,& \\
b_i \,s^2 \,+\, c_i \,t^2 \,+\, d_i\, s t \,+\, X_{\tau(i)}^{[2]} \,+\, A_i\,U^{[2]} \,+\, B_i\, V^{[2]} \,+\,C_i\, W^{[2]} 
&=& 0& (i\le n-2),\\ 
\beta_1 \,s^2 \,+\, \gamma_1\, t^2 \,+\,U^{[2]}&=& 0,& \\
\beta_2 \, s^2 \,+\, \gamma_2\, t^2 \,+\,V^{[2]}&=& 0.&
\end{array}
$$
Thus we obtain
$$
\renewcommand{\arraystretch}{1.0}
\begin{array}{rcccc}
U^{[2]}&=&{\beta_1}\; s^2 &+&{\gamma_1}\; t^2, \\
V^{[2]}&=&{\beta_2}\; s^2 &+& {\gamma_2}\; t^2, \\
W^{[2]}&=&{\alpha_1}\; s^2 &+& {\alpha_2}\; t^2,  
\end{array}
$$
and
$$
\renewcommand{\arraystretch}{1.28}
\begin{array}{ccl}
X_i^{[2]} &=& (b_{\tau(i)}\,s^2 +c_{\tau (i)} t^2 + d_{\tau (i)} \, st ) 
\,+\, A_{\tau (i)} U^{[2]}\,+\, B_{\tau (i)} V^{[2]} \,+\, C_{\tau(i)} W^{[2]}
\\
&=&\phantom{\, +\, }(b_{\tau(i)}\,\,+\, A_{\tau(i)}\,\beta_1 \,+\, B_{\tau(i)}\,\beta_2\,+\, C_{\tau(i)}\,\alpha_1 )\, s^2\\ 
&&\, +\, (c_{\tau(i)}\,\,+\, A_{\tau(i)}\,\rlap{$\gamma_1$}\phantom{\beta_1} \,+
\, B_{\tau(i)}\,\rlap{$\gamma_2$}\phantom{\beta_2}\,+\,C_{\tau(i)}\,\alpha_2)\, t^2\\ &&\,+\, d_{\tau(i)}\,st.
\end{array}
$$
Looking at  the homogeneous part of degree $3$
in $(s, t)$
of $\Phi_L$, we get the equation
$$
\beta_1\, s^3 \,+\, \beta_2\, s^2 t \,+\, \gamma_1\, s t^2 \,+\, \gamma_2\, t^3 \,+\, U^{[2]}\, s\,+\, V^{[2]}\,t \,+\, W^{[3]}\,=\,0.
$$
Thus we obtain
$$
W^{[3]}=0.
$$
Looking at  the homogeneous part of degree $3$ in $(s, t)$
of $(D_{n-1} \Phi)_{L}$ and $(D_n \Phi)_{L}$, we obtain the equations
$$
\begin{array}{lcc}
\dsum _{i=1}^{n-2}d_i X_i^{[2]} t + f_{31} s^2 t + f_{13} t^3 +\;U^{[3]}+
\dsum_{j=n+1}^{m-1} (A_j U^{[2]} +B_j V^{[2]} +C_j W^{[2]})\,e_j t &=&0, \\
\dsum _{i=1}^{n-2}d_i X_i^{[2]} s + f_{31} s^3  + f_{13} s t^2 +V^{[3]}+
\dsum_{j=n+1}^{m-1} (A_j U^{[2]} +B_j V^{[2]} +C_j W^{[2]})\,e_j s &=&0,
\end{array}
$$
where $e_j$ is the coefficient of $\xi_1\xi_2$ in $\phi_j$.
We  put
$$
\renewcommand{\arraystretch}{1.6}
\begin{array}{ccl}
Q(s, t)&:=&\dsum _{i=1}^{n-2}d_i X_i^{[2]}  + f_{31} s^2  + f_{13} t^2 +
\dsum_{j=n+1}^{m-1} (A_j U^{[2]} +B_j V^{[2]} +C_j W^{[2]})\,e_j \\
 &=& Q_s\,s^2\,+\,Q_t\,t^2,
\end{array}
$$
where $Q_s$ and $Q_t$ are given in Table~\ref{table:Q}.
\begin{table}
{
$$
\renewcommand{\arraystretch}{1.0}
\begin{array}{ccl}
Q_s &:=& \dsum_{i=1}^{n-2} (b_{\tau(i)} + A_{\tau(i)}\beta_1 +B_{\tau (i)} \beta_2 + C_{\tau (i)} \alpha_1) d_i+\\
&&\,\,\,+\,\,f_{31}\,+\,\dsum_{j=n+1}^{m-1}(A_j\beta_1+B_j\beta_2+C_j\alpha_1) e_j, \\ 
Q_t &:=& \dsum_{i=1}^{n-2} (c_{\tau(i)} + A_{\tau(i)}\gamma_1 +B_{\tau (i)} \gamma_2 + C_{\tau (i)} \alpha_2) d_i+\\
&&\,\,\,+\,\,f_{13}\,+\,\dsum_{j=n+1}^{m-1}(A_j\gamma_1+B_j\gamma_2+C_j\alpha_2) e_j.
\end{array}
$$
}
\caption{The coefficients $Q_s$ and $Q_t$}\label{table:Q}
\end{table}
(Note that the coefficient of $st$ in $Q(s, t)$ is $\sum_{i=1}^{n-2} d_i d_{\tau (i)}=0$.)
Then  we obtain
$$
U^{[3]}=t\,Q(s, t)\quand V^{[3]}=s\,Q(s, t).
$$
\par
Suppose that $\Lambda$ is defined in $L$ by
\begin{equation*}\label{eq:Lambda}
D u +E v+F w=0,
\end{equation*}
where $D, E, F$ are generic over  $k$.
From the assumption,  $A_i, B_i, C_i$ and $D, E, F$ are algebraically independent over $k$.
The scheme  $\CL$ is defined in $\SL$ by $\Gamma=0$
locally at $P$, where 
$\Gamma$ is the formal power series of $(s, t)$ defined by
$$
\Gamma:=D U +E V+F W.
$$
Obviously, we have $\Gamma^{[0]}=\Gamma^{[1]}=0$.
We also  have
$$
\Gamma^{[2]}\,=\,(D\beta_1 +E\beta_2 +F\alpha_1)\,s^2 \,+\,(D\gamma_1 +E\gamma_2 +F\alpha_2)\, t^2\,=\, \ell^2,
$$
where
$$
\ell\;:=\;\sqrt{D\beta_1 +E\beta_2 +F\alpha_1} \,\, s \; +\;\sqrt{D\gamma_1 +E\gamma_2 +F\alpha_2}\,\, t.
$$
Hence, if $\Gamma\ne0$,
then the curve $\CL$ on $\SL$ has a degenerate singularity at $P$.
Moreover,
$\Gamma^{[2]}$ is not zero
if and only if at least one of 
$\alpha_1, \alpha_2, \beta_1, \beta_2, \gamma_1, \gamma_2$ 
is not zero.
Thus the assertions (1) and (2) are proved.
\par
Suppose that $P\in \fRsm$ and that $\Gamma^{[2]}=\ell^2\ne 0$.
The degenerate singular point $P$ of $\CL$ is an ordinary cusp if and only if
the linear form 
$\ell$ does not divide
$$
\Gamma^{[3]}=(D\, t +E\, s)\, (\sqrt{Q_s} \, s+\sqrt{Q_t}\, t)^2.
$$
The two linear forms $\ell$ and $Dt+Es$ are proportional if and only if
\begin{equation}\label{eq:ell1}
D^2 \, (D\beta_1 +E\beta_2 +F\alpha_1)\;+\;E^2\, (D\gamma_1 +E\gamma_2 +F\alpha_2)\;=\;0
\end{equation}
holds.
Since $D, E, F$ are generic over $k$,
and at least one of $\alpha_1, \alpha_2, \beta_1, \beta_2, \gamma_1, \gamma_2$ 
is not zero by $\Gamma^{[2]}\ne 0$,
the  equality~\eqref{eq:ell1} does not hold.
The two linear forms $\ell$ and $\sqrt{Q_s} s +\sqrt{Q_t} t$ are proportional if and only if
\begin{equation}\label{eq:ell2}
Q_t\, (D\beta_1 +E\beta_2 +F\alpha_1)\;+\;Q_s\, (D\gamma_1 +E\gamma_2 +F\alpha_2)\;\;=\;\;0
\end{equation}
holds.
Since $P\in \fRsm$,
we see from Theorem~\ref{thm:evenmain} (1)
that at least one of
$d_1, \dots, d_{n-2}$, $e_{n+1}, \dots, e_{m-1}$
is not zero.
Expanding the left hand side of~\eqref{eq:ell2}
as a polynomial of $A_i, B_i, C_i$ and $D, E, F$,
we see that~\eqref{eq:ell2} holds for a generic choice of 
$A_i, B_i, C_i,D, E, F$
if and only if the condition~\eqref{eq:rank2}
does \emph{not} hold.
Therefore the assertion (3) is proved.
\end{proof}
\begin{proposition}
Suppose that $\CL$ is of codimension $1$ in $\SL$ at $P$, and is  reduced and irreducible locally at $P$. 
Let $\nu : \wCL\to \CL$ be the normalization of $\CL$ at $P$,
and let $z$  be a formal parameter of $\wCL$ at the point $P\sprime\in \wCL$
such that $\nu (P\sprime)=P$.
Then the image of 
$$
\map{{(\pi_\Lambda\circ\nu)_{P\sprime}\sp *}\comp}{(\OOO_{\Lambda, [f]})\comp}{(\OOO_{\wCL, P\sprime})\comp}
$$
is contained in the subring
$$
T:=\set{f=\sum f_{\mu} z^\mu}{f_1=f_2=f_3=f_5=0}
$$
of $(\OOO_{\wCL, P\sprime})\comp=k[[z]]$.
\end{proposition}
\begin{proof}
We use the same notation as  in the proof of Proposition~\ref{prop:cut}.
The formal completion at $P\sprime$ of the 
composite of $\nu :\wCL\to \CL$ and the inclusion $\iota: \CL\inj \SL$
is given by
$$
\renewcommand{\arraystretch}{1.28}
\begin{array}{ccl}
(\iota\circ \nu)^* s&=& a\, z^2 \,+\, b\, z^3 \,+\,(\hbox{terms of degree $\ge 4$}) \quand \\
(\iota\circ \nu)^* t&=& c\, z^2 \,+\, d\, z^3 \,+\,(\hbox{terms of degree $\ge 4$}),
\end{array}
$$
where $a,b,c,d\in k$.
Since $\Lambda$ is general in $L$, the variables 
$$
u\sprime:= u\rest\Lambda
\quand
v\sprime:= v\rest\Lambda
$$
form a formal parameter system of $\Lambda$ at $[f]$.
The formal completion at $P\sprime$ of the morphism
$\pi_\Lambda\circ \nu: \wCL\to \Lambda$
is given by
$$
\renewcommand{\arraystretch}{1.28}
\begin{array}{ccl}
(\pi_\Lambda\circ \nu)^* u\sprime&=& U(\;a z^2 + b z^3 +\cdots, \;\;c z^2 + d z^3 +\cdots) \quand\\ 
(\pi_\Lambda\circ \nu)^* v\sprime&=& V(\;a z^2 + b z^3 +\cdots, \;\;c z^2 + d z^3 +\cdots).
\end{array}
$$
Since $U^{[1]}$ and $V^{[1]}$ are both zero,
and the coefficients of $st$ in  $U^{[2]}$ and $V^{[2]}$ are also zero,
we see that $(\pi_\Lambda\circ \nu)^* u\sprime$ and $(\pi_\Lambda\circ \nu)^* v\sprime$ are 
contained in the subring $T$ 
of $k[[z]]$.
\end{proof}
\subsection{The case where the linear system is sufficiently ample}
\begin{definition}
For a point $p$ of $\XX$ and a positive integer $d$, we denote by 
$$
\map{v_p\sp{[d]}}{M}{\LLL_p/\mmm_p^{d+1} \LLL_p}
$$
the evaluation homomorphism at $p$,
where $\mmm_p$ is the maximal ideal of  $\OOO_{\XX, p}$,
and $\LLL_p$ is the $ \OOO_{\XX, p}$-module $\LLL\otimes \OOO_{\XX, p}$.
\end{definition}
\begin{proposition}\label{prop:eval}
{\rm (1)}
Suppose that $v_p\sp{[2]}$ is surjective
 at every point $p$ of $\X$.
Then $\fR$ is irreducible.
\par
{\rm (2)}
Suppose that $v_p\sp{[3]}$ is surjective
 at every point $p$ of $\X$.
Then $\fRsm$ is dense in $\fR$.
Moreover,
if $P=\pf$ is a general point of $\fR$,
then the singularity of $D_{[f]}$ at $p$ is of type $\AAeven$ and of type $\typeC$.
\par
{\rm (3)}
Suppose that $v_p\sp{[4]}$ is surjective
 at every point $p$ of $\X$.
If $P=\pf$ is a general point of $\fR$,
then the singularity of $D_{[f]}$ at $p$ is   of type $\typeR$.
\end{proposition}
\begin{proof}
Considering  the first projection $\fR\to\X$,
we deduce  the statements of Proposition from 
Remark~\ref{rem:open}
and Theorem~\ref{thm:evenmain} (1).
\end{proof}
Combining Proposition~\ref{prop:eval}
with the results proved so far,
we obtain the facts (I)-(IV) stated in Introduction.
\section{Divisibility of a Chern number}\label{sec:div}
We assume that $\XXX$ is smooth and 
$\LLL$ is very ample,
and put $|M|=|\LLL|$.
Then we have $\X=\XX=\XXX$ and $\fCC=\fC$.
We further assume that $\fRsm$ is dense in $\fR$.
Under these assumptions, we calculate
$$
\deg \fR:= \int_{\XXP} c_1 (\pr_2^* \OOO_{\PM} (1) )^{m-2} \cap [\fR].
$$
We put
$$
\lambda:=c_1 (\pr_1^* \LLL)\quand
h:= c_1 (\pr_2^* \OOO_{\PM} (1) ).
$$
Then we have
$$
c_1 (\widetilde \LLL)=h+\lambda.
$$
For simplicity, we write $c_i (X)$ for $\pr_1^* c_i (X)$.
By Proposition~\ref{prop:fCC} and Thom-Porteous formula~\cite[Chapter 14]{MR1644323},
we have
\begin{equation*}\label{eq:fCC}
[\fCC]= 
\left((\lambda + h)\, \sum\sb{i=0}\sp{n} (-1)^i c_i(X) (\lambda +h)^{n-i} \right) \,\cap\, [\XXP]
\end{equation*}
in the Chow group of $\XXP$.
By Theorem~\ref{thm:HTP},
we have
$$
[\fR]= \left(\frac{n}{2}\,(h+\lambda)  - c_1 (X)\right) \,\cap\,[\fCC].
$$
We put
$$
G:=\left(\frac{n}{2}\,(h+\lambda)  - c_1 (X)\right) \left((\lambda + h)\, \sum\sb{i=0}\sp{n} (-1)^i c_i(X) (\lambda +h)^{n-i} \right).
$$
Then   we obtain 
$$
\deg \fR=\int_{\XX} (\hbox{the coefficient of $h^2$ in  $G$})\,\cap\,[X].
$$
By Corollary~\ref{cor:factorKKK},
the integer
$\deg \fR$
is divisible by $4$.
If the coefficient of $h^2$ in  $G$ is regarded as a polynomial of $\lambda$,
then the constant term is equal to  $ (-1)^n (n\,c_n(X) /2 + c_1 (X) c_{n-1} (X))$.
Putting $\LLL=\AAA\sp{\otimes {4}}$
with $\AAA$ a sufficiently ample line bundle, 
we obtain the following:
\begin{corollary}
Let $X$ be a smooth projective variety in characteristic $2$ of dimension $n$ being even.
Then the integer
$$
\int_X \left(\frac{n}{2}\, c_n(X) + c_1 (X) c_{n-1} (X)\right)\cap [X]
$$
is divisible by $4$.
\end{corollary}
In fact, this divisibility relation follows  from the Hirzebruch-Riemann-Roch theorem
by the argument of Libgober and Wood.
See~\cite[Remark 2.4]{MR1064869}.

\def\cprime{$'$} \def\cprime{$'$}
\providecommand{\bysame}{\leavevmode\hbox to3em{\hrulefill}\thinspace}

\end{document}